\documentstyle[amsmath,amssymb,11pt]{article}

\newcommand{\rom}[1]{{\rm #1}}
\newcommand{\wick}[1]{{:}\omega^{\otimes #1}{:}_\lambda}

\allowdisplaybreaks

\makeatletter\@addtoreset{equation}{section}\makeatother

\begin{document}

\setcounter{page}{1} \setcounter{section}{0} \thispagestyle{empty}

\newtheorem{definition}{Definition}[section]
\newtheorem{remark}{Remark}[section]
\newtheorem{proposition}{Proposition}[section]
\newtheorem{theorem}{Theorem}[section]
\newtheorem{corollary}{Corollary}[section]
\newtheorem{lemma}{Lemma}[section]

\newcommand{\indlim}{\operatornamewithlimits{ind\,lim}}
\newcommand{\Ffin}{{\cal F}_{\mathrm fin}}
\newcommand{\Fext}{{\cal F}_{\mathrm ext}}
\newcommand{\D}{{\cal D}}
\newcommand{\N}{{\Bbb N}}
\newcommand{\C}{{\Bbb C}}
\newcommand{\Z}{{\Bbb Z}}
\newcommand{\R}{{\Bbb R}}
\newcommand{\Rp}{{\R_+}}
\newcommand{\eps}{\varepsilon}
\newcommand{\supp}{\operatorname{supp}}
\newcommand{\la}{\langle}
\newcommand{\ra}{\rangle}
\newcommand{\const}{\operatorname{const}}
\renewcommand{\emptyset}{\varnothing}
\newcommand{\di}{\partial}
\newcommand{\hotimes}{\hat\otimes}

\renewcommand{\author}[1]{\medskip{\Large #1}\par\medskip}

\begin{center}{\Large \bf Polynomials of Meixner's type  in
infinite dimensions---\\[2mm]Jacobi fields and orthogonality
measures}\end{center}

\author{Eugene Lytvynov}

\noindent{\sl Institut f\"{u}r Angewandte Mathematik,
Universit\"{a}t Bonn, Wegelerstr.~6, D-53115 Bonn,\\ Germany; SFB
611, Univ.~Bonn, Germany;  BiBoS, Univ.\ Bielefeld, Germany\\[2mm]
 {\rm E-mail: lytvynov@wiener.iam.uni-bonn.de}}

\begin{abstract}
\noindent The classical polynomials of Meixner's type---Hermite,
Charlier, Laguerre, Meixner, and Meixner--Pollaczek
polynomials---are distinguished through a special form of their
generating function, which involves the Laplace transform of their
orthogonality measure. In this paper, we study analogs of the
latter three classes of polynomials in infinite dimensions. We fix
as an underlying space a (non-compact) Riemannian manifold $X$ and
an intensity measure $\sigma$ on it. We consider a Jacobi field in
the extended Fock space over $L^2(X;\sigma)$, whose field operator
at a point $x\in X$ is of the form $\di_x^\dag+\lambda\di_x^\dag
\di_x+\di_x+\di^\dag_x\di_x\di_x$, where $\lambda$ is a real
parameter. Here, $\di_x$ and $\di_x^\dag$  are, respectively, the
annihilation and creation operators at the point $x$. We then
realize the field operators as multiplication operators in
$L^2({\cal D}';\mu_\lambda)$, where ${\cal D}'$ is the dual of
${\cal D}{:=}C_0^\infty(X)$, and $\mu_\lambda$ is the spectral
measure of the Jacobi field. We show that $\mu_\lambda$ is a gamma
measure  for $|\lambda|=2$, a Pascal measure for $|\lambda|>2$,
and a   Meixner measure for $|\lambda|<2$. In all the cases,
$\mu_\lambda$  is a L\'evy  noise measure. The isomorphism between
the extended Fock space and $L^2({\cal D}';\mu_\lambda)$ is
carried out by infinite-dimensional polynomials of Meixner's type.
We find the generating function of these polynomials and using it,
we study the action of the operators $\di_x$ and $\di_x^\dag$ in
the functional realization.
\end{abstract}


\section{Introduction}

In his classical work \cite{Meixner}, Meixner considered the
following problem: Suppose that functions $f(z)$ and $\Psi(z)$ can
be expanded in a formal power series of $z\in\C$ and suppose that
$f(0)=1$, $\Psi(0)=0$, and $\Psi'(0)=1$. Then, the equation
\begin{equation}\label{jxaboiu}
G(x,z){:=}\exp(x\Psi(z))f(z)=\sum_{n=0}^\infty
\frac{P^{(n)}(x)}{n!}\,z^n
\end{equation} generates a system of polynomials $P^{(n)}(x)$,
$n\in\Z_+$, with leading coefficient 1.  (These polynomials are
now called Sheffer polynomials.) Find all polynomials of such type
which are orthogonal with respect to some probability measure
$\mu$
 on $\R$. To solve this problem, Meixner essentially used the two
 following facts. First, by the Favard theorem, a system of
 polynomials $P^{(n)}(x)$, $n\in\Z_+$, with leading coefficient
 1 is  orthogonal if and only if these polynomials satisfy the
 recurrence formula
\begin{equation}\label{hft89uioj}xP^{(n)}(x)=P^{(n+1)}(x)+a_nP^{(n)}(x)+b_n
P^{(n-1)}(x),\qquad n\in\Z_+,\ P^{(-1)}(x){:=}0,\end{equation}
with real numbers $a_n$ and positive numbers $b_n$; or
equivalently, the polynomials $P^{(n)}(x)$ determine the infinite
Jacobi matrix with the elements $a_n$ on the main diagonal and the
elements $\sqrt{b_n}$ on the upper and lower diagonals. And
second, as follows from \eqref{jxaboiu},
\begin{equation}\label{jkrsrh}
\Psi^{-1}(D)P^{(n)}(x)=nP^{(n-1)}(x),\qquad n\in\N,\end{equation}
where $\Psi^{-1}$ is the inverse function of $\Psi$ and
$D{:=}\frac d{dx}$. Meixner showed  that the solution of this
problem is completely determined by the equations $$
\lambda=a_n-a_{n-1},\quad n\in\N,\qquad \kappa
=\frac{b_n}n-\frac{b_{n-1}}{n-1},\quad n\ge2,$$ where $\lambda$
and $\kappa$ are some parameters. If $\kappa=0$, we have to
distinguish the two  following cases:

 I) $\lambda=0$; without loss of generality, we then get  $a_n=0$
and $b_n=n$ in \eqref{hft89uioj}, $P^{(n)}(x)$  are the Hermite
polynomials, $\mu $ is the standard Gaussian distribution on $\R$.

 II) $\lambda\ne0$, so that
$a_n=\lambda n$, $b_n=n$,  $P^{(n)}(x)$  are the Charlier
polynomials, $\mu $ is   a centered Poisson distribution on $\R$.

Let now $\kappa\ne0$ and we set $\kappa=1$ for simplicity of
notations. We then get $a_n=\lambda n$ and $b_n=n^2$. We introduce
two quantities $\alpha$ and $\beta$ through the equation
\begin{equation}\label{jkxhgvfuz} 1+\lambda z+ z^2=(1-\alpha
z)(1-\beta z).\end{equation} We now have to distinguish the three
following cases:

III) $|\lambda|=2$, so that $\alpha=\beta=\pm 1$, $P^{(n)}(x)$ are
the Laguerre polynomials, $\mu$ is a centered gamma distribution.

IV) $|\lambda|>2$, so that $\alpha\ne\beta$, both real,
$P^{(n)}(x)$ are the Meixner polynomials (of the first kind),
which are orthogonal with respect to a centered Pascal (negative
binomial) distribution.

V) $|\lambda|<2$, so that $\alpha\ne\beta$, both complex
conjugate,   $P^{(n)}(x)$ are the Meixner polynomials of the
second kind, or the Meixner--Pollaczek polynomials in other terms.
These  are orthogonal with respect to a measure $\mu$ obtained by
centering  a probability measure  of the form
$C\exp(ax)|\Gamma(1+im x )|^2\,dx$, where $a\in\R$, $m>0$, and $C$
is the normalizing constant. We will call it a Meixner measure,
though there seems to be no commonly accepted name for it.

In all the above cases, the generating function $G(x,z)$ defined
in \eqref{jxaboiu} can be represented as
$G(x,z)=\exp\big(x\Psi(z)\big)/L_\mu(\Psi(z))$, where
$L_\mu(z){:=}\int_\R e^{zx}\,\mu(dx)$ is the extension of the
Laplace transform of the measure $\mu$ defined in a neighborhood
of zero in $\C$.

In  the present paper, we will  study  analogs of polynomials of
Meixner's type and their orthogonality measures in infinite
dimensions. In the case of the Gaussian and Poisson measures, such
a theory  is, of course, well studied; we refer to e.g. \cite{BK,
HKPS} for the Gaussian case and to  e.g.\ \cite{KKO,KSS} for the
Poisson case. Notice that, in both cases, the Fock space and the
corresponding Jacobi fields of operators in it play a fundamental
role (see \cite{bere,BeLi,Ly} for the notion of a Jacobi field in
the Fock space). In particular, the field operator at a point
$x\in X$, where $X$ is an underlying space, has the form
$\di_x^\dag+\di_x$ in the Gaussian case, and
$\di_x^\dag+\lambda\di_x^\dag\di_x+\di_x$ in the Poisson case.
Here, $\di_x$ and $\di_x^\dag$ are  the annihilation and  creation
operators at the point $x$, respectively.

Concerning the gamma case, III), an infinite-dimensional analog of
the Laguerre polynomials and the corresponding Jacobi field was
studied in \cite{KL, silva}. The polynomials are now orthogonal
with respect to  the (infinite-dimensional) gamma measure, which
is a special case of a compound Poisson measure. Since such a
measure does not possess the chaotic decomposition property,
instead of the usual Fock space one has to use the so-called
extended Fock space. This space,   on the one hand, extends the
usual Fock space and, on the other hand, still has some
similarities with it. The field operator at a point $x\in X$ has
the form $\di_x^\dag+2\di^\dag_x\di_x+\di_x+\di^\dag_x\di_x\di_x$.
In \cite{beme1}, the structure of the extended Fock space was
discussed in detail, and  in \cite{beme2}, it was shown that the
extended Fock space decomposition of the Gamma process can be
thought  of as an expansion of any $L^2$-random variable in
multiple integrals constructed by using a family of resolutions of
the identity in the extended Fock space. We also refer to the
recent paper \cite{Tsil} and the references therein, where many
other properties of the gamma measure are discussed in detail.

As for the cases IV) and V), the role of the orthogonality measure
should be played by (infinite-dimensional) Pascal and Meixner
measures (processes). These processes in the case $X=\R_+$, both
L\'evy, were introduced in \cite{BR} and \cite{ST}, respectively.
In \cite{Gri}, the Meixner process was proposed for a model for
risky asserts and an analog of the Black--Sholes formula was
established. In \cite{NS} (see also the recent book \cite{S}), the
gamma, Pascal, and Meixner processes served as main examples of a
chaotic representation for every square-integrable random variable
in terms of the  orthogonalized Teugels martingales related to the
process. Though the one-dimensional polynomials of Meixner's type
were used in this work in order to carry out the orthogonalization
procedure of the Teugels martingales (which, in turn,  are the
centered power jump processes related to the original process),
infinite-dimensional polynomials corresponding to these processes
have not appeared in this work; furthermore, they were mentioned
as an open problem in \cite{S}.

The contents of the present paper is as follows. In
Section~\ref{uisdygfuzuzg}, we fix as an underlying space $X$  a
smooth (non-compact) Riemannian manifold and an intensity measure
$\sigma$ on it. We consider a Jacobi field in the extended Fock
space over $L^2(X;\sigma)$, whose field operator at a point $x\in
X$ has the form
$\di_x^\dag+\lambda\di^\dag_x\di_x+\di_x+\di^\dag_x\di_x\di_x$,
where $\lambda$ is a fixed real parameter. Using ideas of
\cite{BeLi,KL,Ly}, we construct via the projection spectral
theorem \cite{BK} a Fourier transform in generalized joint
eigenvectors of the Jacobi field. This gives us a unitary operator
$I_\lambda$ between the extended Fock space and the space
$L^2({\cal D}';\mu_\lambda)$, where ${\cal D}'$ is the dual space
of ${\cal D}{:=}C_0^\infty(X)$ with respect to the zero space
$L^2(X;\sigma)$, and $\mu_\lambda$ is the spectral measure of the
Jacobi field, i.e., the image of any field operator  under
$I_\lambda$ is a multiplication operator in $L^2({\cal
D}';\mu_\lambda)$. The $\mu_\lambda$ is a gamma measure  for
$|\lambda|=2$, a Pascal measure for $|\lambda|>2$, and a Meixner
measure for $|\lambda|<2$,  in the sense that, for any bounded
$\Delta\subset X $, the (naturally defined) random variable $\la
\cdot,\chi_\Delta\ra$ has a corresponding one-dimensional
distribution. Furthermore, for $|\lambda|\ge2$ $\mu_\lambda$ is a
compound Poisson measure, and for $|\lambda|<2$ $\mu_\lambda$ is a
L\'evy noise measure. In particular, for $X=\R$ we obtain the
gamma, Pascal, and Meixner processes, respectively.

Next, under the unitary $I_\lambda$, the image of any vector
$f^{(i)}\in{\cal D}_\C^{\hotimes i}$ is a continuous polynomial
$\la\wick i,f^{(i)}\ra$ of the variable $\omega\in{\cal D}'$,
which may be understood as an infinite-dimensional polynomial of
Meixner's type, since $\la\wick i,\chi_\Delta^{\otimes
i}\ra=P^{(i)}_{\lambda,\,\Delta}(\la\omega,\chi_\Delta\ra)$, where
$P^{(i)}_{\lambda,\,\Delta}(\cdot)$ is a one-dimensional
polynomial  of Meixner's type.

In Section~\ref{jfrtd}, we identify the generating function
$G_\lambda(\omega,\varphi){:=}\sum_{n=0}^\infty
\frac1{n!}\,\la\wick n,\varphi^{\otimes n}\ra$, $\omega\in{\cal
D}'$, $\varphi\in{\cal D }_\C$, and show that
$G_\lambda(\omega,\varphi)=\exp\big(\la\omega,\Psi_\lambda(\varphi)\ra\big)/
\ell_\lambda(\Psi_\lambda(\varphi))$, where $\ell_\lambda$ is the
extension of the Laplace transform of the measure $\mu_\lambda$
defined in a neighborhood of zero in ${\cal D}_\C$, and
$\Psi_\lambda$ is the same function as $\Psi$ in \eqref{jxaboiu}.

Finally, in Section~\ref{gtzdfzrdr}, using results of
\cite{KSS,KSWY}, we introduce a test space $(\cal D_\lambda)$
consisting of those functions on ${\cal D}'$ which may be extended
to  entire functions on ${\cal D}_\C'$ of first order of growth
and of minimal type. This space is densely and continuously
embedded into $L^2({\cal D}';\mu_\lambda)$.  We then study the
action of the operators $\di_x:({\cal D}_\lambda)\to({\cal
D}_\lambda)$ and $\di_x^\dag:({\cal D}_\lambda)\to ({\cal
D}_\lambda)^*$, where $({\cal D}_\lambda)^*$ is the dual of
$({\cal D_\lambda})$. We note that, analogously to \eqref{jkrsrh},
we have $\di_x=\Psi_\lambda^{-1}(\nabla_x)$, where $\nabla_x$ is
the G\^ateaux derivative in direction $\delta_x$. We obtain
explicit formulas for the operators $\di_x$ and $\int_X
\sigma(dx)\,\xi(x)\di^\dag_x$, $\xi\in{\cal D}$. It should be
stressed that, for the latter operator in the case $|\lambda|<2$,
the possibility of a (unique) extension of a test function on
${\cal D}'$  to a function on ${\cal D}_\C'$ plays a principle
role.

In a forthcoming paper, we will study a connection between the
extended Fock space decomposition of $L^2({\cal D}';\mu_\lambda)$
obtained in this paper and the chaotic decomposition of this space
in the case $X=\R$ as in \cite{NS}. Finally, we note that one can
also study a more general model where, in the field operator at a
point $x\in X$, the value of the parameter $\lambda$ depends on
$x$. Then, the corresponding noise will be with independent values
and at each point $x\in X$ its properties will be the same as the
properties of the  noise at the point $x$ under
$\mu_{\lambda(x)}$.

\section{Meixner's Jacobi field and its spectral
measures}\label{uisdygfuzuzg}

Let $X$ be a complete, connected, oriented $C^\infty$
(non-compact) Riemannian manifold and let ${\cal B}(X)$ be the
Borel $\sigma$-algebra on $X$. Let $\sigma$ be a Radon measure on
$(X,{\cal B}(X))$ that is non-atomic, i.e., $\sigma(\{x\})=0$ for
every $x\in X$ and non-degenerate, i.e., $\sigma(O)>0$ for any
open set $O\subset X$. (We note the the assumption of the
nondegeneracy of $\sigma$ is nonessential and the results below
may be generalized to the case of a degenerate $\sigma$.) Note
that $\sigma(\Lambda)<\infty$ for each $\Lambda\in{\cal
O}_c(X)$---the set of all open  sets in $X$ with compact closure.

We denote by ${\cal D}$ the space $C_0^\infty(X)$ of all
real-valued infinite differentiable functions on $X$ with compact
support. This space may be naturally endowed with a topology of a
nuclear space, see e.g.\ \cite{BUS} for the case $X=\R^d$ and
e.g.\  \cite{Die} for the case of a general Riemannian manifold.
We recall that
\begin{equation}\label{uiguz}{\cal D}=\operatornamewithlimits{proj\,lim}_{\tau\in
T}{\cal H}_\tau.\end{equation} Here, $T$ denotes the set of all
pairs $(\tau_1,\tau_2)$ with $\tau_1\in\Z_+$ and $\tau_2\in
C^\infty(X)$, $\tau_2(x)\ge1$ for all $x\in X$, and ${\cal
H}_{\tau}={\cal H}_{(\tau_1,\tau_2)}$ is the Sobolev space on $X$
of order $\tau_1$ weighted by the function $\tau_2$, i.e., the
scalar product in ${\cal H}_\tau$, denoted by $(\cdot,\cdot)_\tau$
is given by \begin{equation}\label{hjdrt}(f,g)_{\tau}=\int_X
\bigg( f(x)g(x)+\sum_{i=1}^{\tau_1} \langle
\nabla^{i}f(x),\nabla^{i}g(x)\ra_{T_x(X)^{\otimes{i}}}\bigg)\tau_2(x)\,dx,\end{equation}
where $\nabla^i$ denotes the $i$-th (covariant) gradient, and $dx$
is the volume measure on $X$. For $\tau,\tau'\in T$, we will write
$\tau'\ge\tau$ if $\tau'_1\ge\tau_1$ and  $\tau'_2(x)\ge\tau_2(x)$
for all $x\in X$.

 The space $\cal D$ is densely and continuously
embedded into $L^2(X;\sigma)$. As easily seen, there always exists
$\tau_0\in T$ such that ${\cal H}_{\tau_0}$ is continuously
embedded into $L^2(X;\sigma)$. We denote $T'{:=}\{\tau\in T:
\tau\ge \tau_0\}$ and  \eqref{uiguz} holds with $T$ replaced by
$T'$. In what follows, we will just write $T$ instead of $T'$. Let
${\cal H}_{-\tau}$ denote the dual space of ${\cal H}_\tau$ with
respect to the zero space ${\cal H}{:=}L^2(X;\sigma)$. Then ${\cal
D}'=\operatornamewithlimits{ind\,lim}_{\tau\in T}{\cal H}_{-\tau}$
is the dual of ${\cal D}$ with respect to ${\cal H}$, and we thus
get the standard triple $$ {\cal D}'\supset {\cal H} \supset {\cal
D}.$$ The dual pairing between any $\omega\in{\cal D}'$ and
$\xi\in{\cal D}$ will be denoted by $\la\omega,\xi\ra$.

Following \cite{KL}, we define, for each $n\in\N$, an $n$-particle
extended Fock space over ${\cal H}$, denoted by $\Fext^{(n)}({\cal
H})$. Under a loop $\kappa$ connecting points $x_1,\dots,x_m$,
$m\ge2$, we understand a class of ordered sets
$(x_{\pi(1)},\dots,x_{\pi(m)})$, where $\pi$ is a permutation of
$\{1,\dots,n\}$, which coincide up to a cyclic permutation. We put
$|\kappa|=m$. We will also interpret a set $\{x\}$ as a
``one-point'' loop $\kappa$, i.e., a loop that comes out of $x$,
$|\kappa|=1$. Let
$\alpha_n=\{\kappa_1,\dots,\kappa_{|\alpha_n|}\}$ be a collection
of loops $\kappa_j$ that connect points from the set
$\{x_1,\dots,x_n\}$ so that every point $x_i\in\{x_1,\dots,x_n\}$
goes into one loop $\kappa_j=\kappa_{j(i)}$ from $\alpha_n$. Here,
$|\alpha_n|$ denotes the number of the loops in $\alpha_n$,
evidently $n=\sum_{j=1}^{|\alpha_n|}|\kappa_j|$. Let $A_n$ stand
for the set of all possible collections of loops $\alpha_n$ over
the points $\{x_1,\dots,x_n\}$. (We note that the set $A_n$
contains $n!$ elements  \cite[Remark~2.1]{KL}.)  Every
$\alpha_n\in A_n$ generates the following continuous mapping
\begin{equation}\label{1.2}
 {\cal D}^{\hotimes  n}_\C\ni f^{(n)}=f^{(n)}(x_1,\dots,x_n)
\mapsto f^{(n)}_{\alpha_n}(
\underbrace{x_1,\dots,x_1}_{\text{$|\kappa_1|$ times}} ,
\underbrace{x_2,\dots,x_2}_{\text{$|\kappa_2|$ times}} ,\dots,
\underbrace{x_{|\alpha_n|},
\dots,x_{|\alpha_n|}}_{\text{$|\kappa_{|\alpha_n|}|$ times}} )\in
{\cal D}_\C^{\otimes|\alpha_n|},
           \end{equation}
where the lower index $\C$ denotes complexification of a real
space and the symbol $\hotimes$ stands for the symmetric tensor
power. We define a scalar product on ${\cal D}_{\C}^{\hotimes n}$
by
\begin{equation}\label{1.5}
( f^{(n)},g^{(n)})_{\Fext^{(n)}({\cal H})}=\sum_{\alpha_n\in
A_n}\int _{X^{|\alpha_n|}}\big(\overline{
f^{(n)}}g^{(n)}\big)_{\alpha_n}\, d\sigma^{\otimes |\alpha_n|},
\end{equation}
where $\overline{f^{(n)}}$ is the complex conjugate of $f^{(n)}$.
 Let $\Fext^{(n)}({\cal
H})$ be the closure of ${\cal D}_{\C}^{\hotimes n}$ in the norm
generated by \eqref{1.5}.

The extended Fock space $ \Fext({\cal H})$ over ${\cal H}$ is
defined as a weighted direct sum of the spaces $\Fext^{(n)}({\cal
H })$:
\begin{equation}\label{1.5a}
\Fext({\cal H})=\bigoplus_{n=0}^\infty \Fext^{(n)}({\cal H
})\,n!,\end{equation} where $\Fext^{(0)}({\cal H})=\C$ and $0!=1$.
That is, $\Fext({\cal H})$ consists of sequences $f=(f^{(0)},
f^{(1)},$ $f^{(2)},\dots)$ such that $f^{(n)}\in \Fext^{(n)}({\cal
H })$ and $$ \|f\|_{\Fext({\cal H
})}^2=\sum_{n=0}^\infty\|f^{(n)}\|^2_{\Fext^{(n)}({\cal H})}
n!<\infty.$$ We will always identify any
$f^{(n)}\in\Fext^{(n)}({\cal H})$ with the element
$(0,\dots,0,f^{(n)},0,\dots)$ of $\Fext({\cal H})$.

Note that the usual Fock space ${\cal F}({\cal H})$ can be
considered as a subspace of $\Fext({\cal H})$ generated by
functions $f^{(n)}\in{\cal D}_\C^{\hotimes n}$ such that
$f^{(n)}(x_1,\dots,x_n)=0$ if $x_i=x_j$ for some
$i,j\in\{1,\dots,n\}$, $i\ne j$.

Let $\Ffin(\cal D)$ denote the topological direct sum of the
spaces ${\cal D}_\C^{\hotimes n}$, i.e., $\Ffin ({\cal D})$
consists of all sequences
$f=(f^{(0)},f^{(1)},\dots,f^{(m)},0,0,\dots)$ such that
$f^{(n)}\in{\cal D}_\C^{\hotimes n}$ and the convergence in
$\Ffin({\cal D})$ means uniform finiteness and the coordinate-wise
convergence. Since each space ${\cal D}_\C^{\hotimes n}$ is
nuclear, so is $\Ffin({\cal D})$. As easily seen, the space
$\Ffin(\cal D)$ is densely and continuously embedded into
$\Fext({\cal H})$.

For each $\xi\in{\cal D}$, let $a^+(\xi)$ be the standard creation
operator defined on $\Ffin({\cal D})$:
$$a^+(\xi)f^{(n)}=\xi\hotimes f^{(n)}, \qquad f^{(n)}\in{\cal
D}_\C^{\hotimes n}, \ n\in\Z_+.$$ A simple calculation shows that
the adjoint operator of $a^+(\xi)$ in $\Fext({\cal H})$,
restricted to $\Ffin({\cal D})$, is given by the formula $$
a^-(\xi)=(a^+(\xi))^*\restriction\Ffin({\cal
D})=a_1^-(\xi)+a_2^-(\xi),$$ where $a_1^-(\xi)$ is the standard
annihilation operator: $$
(a_1^-(\xi)f^{(n)})(x_1,\dots,x_{n-1})=n\int_X\xi(x)f^{(n)}(x,x_1,\dots,x_{n-1})\,\sigma(dx),$$
and $a_2^-(\xi)$ is given by
$$(a_2^-(\xi)f^{(n)})(x_1,\dots,x_{n-1})=n(n-1)\big(\xi(x_1)f^{(n)}(x_1,x_1,x_2,\dots,x_{n-1})\big)^\sim,$$
where $(\cdot)^\sim$ denotes the symmetrization of a function.

Finally, we define on $\Ffin({\cal D})$ the neutral operator
$a^0(\xi)$, $\xi\in{\cal D}$, in a standard way: $$
(a^0(\xi)f^{(n)})(x_1,\dots,x_n)=n\big(\xi(x_1)f^{(n)}(x_1,\dots,x_n)\big)^\sim.$$
One easily checks that $a^0(\xi)$ is a Hermitian operator in
$\Fext({\cal H})$.

Now, we fix a parameter $\lambda\in[0,\infty)$ and define
operators $$ a_\lambda(\xi){:=}a^+(\xi)+\lambda a^0(\xi) +
a^-(\xi)+ c_\lambda \la\xi\ra\operatorname{id},\qquad \xi\in{\cal
D}.$$ Here, $\la\xi\ra{:=}\int_X\xi(x)\,\sigma(dx)$,
$\operatorname{id}$ denotes the identity operator, and the
constant $c_\lambda>0$ is given by
\begin{equation}\label{qwqwuzzuuz} c_\lambda{:=}\begin{cases}
\lambda /2,&\text{if } \lambda\in[0,2],\\
2/(\lambda+\sqrt{\lambda^2-4}),&\text{if
}\lambda>2.\end{cases}\end{equation} (The special choice of this
constant will become clear later on, however it is not of a  real
importance.) Each $a_\lambda(\xi)$ with domain $\Ffin({\cal D})$
is a Hermitian operator in $\Fext({\cal H})$.

By construction, the family  of operators
$(a_\lambda(\xi))_{\xi\in{\cal D}}$ has a Jacobi field structure
in the extended Fock space $\Fext({\cal H})$ (compare with
\cite{bere,BeLi,Ly}). We will call this family Meixner's Jacobi
field corresponding to the parameter $\lambda$.

\begin{lemma} The operators $a_\lambda(\xi)$\rom, $\xi\in{\cal
D}$\rom, with domain $\Ffin({\cal D})$ are essentially selfadjoint
in $\Fext({\cal H})$\rom, and their closures $a_\lambda^\sim(\xi)$
constitute a family of commuting selfadjoint operators\rom, where
the commutation is understood in the sense of the commutation of
their resolutions of the identity\rom.

\end{lemma}

\noindent{\it Proof}. The lemma follows directly from
\cite[Theorem~4.1]{bere} whose proof admits a direct
generalization to the case of the extended Fock space.\quad
$\blacksquare$

\begin{theorem}\label{z8fzvear} For each $\lambda\in[0,\infty)$, there exist a
unique probability measure $\mu_\lambda$ on $({\cal D}',{\cal
C}_\sigma({\cal D}'))$\rom, where ${\cal C}_\sigma({\cal D}')$ is
the cylinder $\sigma$-algebra on ${\cal D}'$\rom, and a unique
unitary operator $$I_\lambda:\Fext({\cal H})\to L^2({\cal
D}';\mu_\lambda)$$ such that, for each $\xi\in{\cal D}$\rom, the
image of $a_\lambda^\sim(\xi)$ under $I_\lambda$ is the operator
of multiplication by the function $\la\cdot,\xi\ra$ in $L^2({\cal
D}';\mu_\lambda)$ and $I_\lambda\Omega=1$\rom, where
$\Omega{:=}(1,0,0,\dots)$\rom. The Fourier transform of the
measure $\mu_\lambda$ is given, in a neighborhood of zero, by the
following formula\rom: for  $\lambda=2$
\begin{equation}\label{jasuetfzcw}\int_{{\cal
D}'}e^{i\la\omega,\varphi\ra}\,d\mu_2(\omega)=\exp\left[-\int_X
\log(1-i\varphi(x))\, \sigma(dx)\right],\qquad \varphi\in{\cal
D},\ \|\varphi\|_\infty{:=}\sup_{x\in
X}|\varphi(x)|<1,\end{equation} and for $\lambda\ne2$
\begin{equation}\label{hisav} \int_{{\cal
D}'}e^{i\la\omega,\varphi\ra}\,d\mu_\lambda(\omega)=\exp\left[-\frac1{\alpha\beta}\int_X\log\bigg(
\frac{\alpha e^{-i\beta\varphi(x)}-\beta
e^{-i\alpha\varphi(x)}}{\alpha-\beta}\bigg)\,\sigma(dx)+ic_\lambda\la\varphi\ra\right]
\end{equation} for all $\varphi\in{\cal D}$ satisfying \begin{equation} \label{dyuivf}\bigg\| \frac{\alpha
(e^{-i\beta\varphi}-1)-\beta(e^{-i\alpha\varphi}-1)
}{\alpha-\beta}\bigg\|_\infty<1,\end{equation} $\alpha,\beta$
defined by \eqref{jkxhgvfuz}\rom. The unitary operator $I_\lambda$
is given on the dense set $\Ffin({\cal D})$ by the formula
$$\Ffin({\cal D})\ni f=(f^{(n)})_{n=0}^\infty \mapsto I_\lambda
f=(I_\lambda f)(\omega)=\sum_{n=0}^\infty \la \wick n,f^{(n)}\ra$$
\rom(the series is\rom, in fact\rom, finite\rom)\rom, where $\wick
n\in{\cal D}^{'\hotimes n}$ is defined by the recurrence formula
\begin{gather}\wick{(n+1)}=\wick{n+1}(x_1,\dots,x_{n+1})=\big(\wick
n(x_1,\dots,x_n)\omega(x_{n+1})\big)^\sim\notag\\ \text{}
-n\big(\wick{(n-1)}(x_1,\dots,x_{n-1})\delta(x_{n+1}-x_n)\big)^\sim
\notag\\ \text{} -n(n-1)\big(\wick
{(n-1)}(x_1,\dots,x_{n-1})\delta(x_n-x_{n-1})\delta(x_{n+1}-x_n)\big)^\sim\notag\\
\text{} -\lambda n\big(\wick n
(x_1,\dots,x_n)\delta(x_{n+1}-x_n)\big)^\sim-c_\lambda \big(\wick
n(x_1,\dots,x_n)1(x_{n+1})\big)^\sim,  \notag\\ \wick 0 =1,\ \wick
1=\omega-c_\lambda.\label{udsH} \end{gather}
\end{theorem}

\begin{remark}
\rom{It can be shown that  $\mu_\lambda$ is the spectral measure
of the commutative family of selfadjoint operators
$a^\sim_\lambda(\xi)$, $\xi\in{\cal D}$ (see \cite[Ch.~3]{BK} for
the notion of a spectral measure).}
\end{remark}

\begin{remark}
\rom{Note that taking a parameter $\lambda<0$ would lead us to the
measure $\mu_\lambda$ obtained from the measure $\mu_{-\lambda}$
by the  transformation $\omega\mapsto-\omega$ of the space ${\cal
D}'$, which is why we have excluded this choice.}
\end{remark}

\noindent {\it Proof of Theorem~}\ref{z8fzvear}. As easily seen,
for any $\xi\in{\cal D}$ and $n\in\N$, the operators $a^+(\xi)$,
$a^0(\xi)$, and $a^-(\xi)$ act continuously from ${\cal
D}_\C^{\hotimes n}$ into
 ${\cal D}_\C^{\hotimes (n+1)}$,  ${\cal D}_\C^{\hotimes n}$,  ${\cal D}_\C^{\hotimes
 (n-1)}$, respectively. Therefore, for any $\xi\in{\cal D}$, $a_\lambda(\xi)$
acts continuously on $\Ffin({\cal D})$. Furthermore, for any fixed
$f\in\Ffin({\cal D})$, the mapping ${\cal D}\ni\xi\mapsto
a_\lambda(\xi)f\in\Ffin({\cal D})$ is linear and continuous.
Finally, the vacuum $\Omega$ is evidently a cyclic vector for the
operators $a^+(\xi)$, $\xi\in{\cal D}$, in $\Ffin({\cal D})$, and
hence in $\Fext({\cal H})$. Then, using the Jacobi filed structure
of $a_\lambda(\xi)$, it is easy to show that $\Omega$ is a cyclic
vector for $a_\lambda(\xi)$, $\xi\in{\cal D}$, in $\Fext({\cal
H})$. Thus, analogously to  \cite{BeLi,Ly}, and
\cite[Theorem~3.1]{KL}, we deduce, by using the projection
spectral theorem \cite[Ch.~3, Th.~2.7 and subsec.~3.3.1]{BK}, the
existence of a unique probability measure $\mu_\lambda$ on $({\cal
D}',{\cal C}_\sigma({\cal D}'))$ and a unique unitary operator
$I_\lambda:\Fext({\cal H})\to L^2({\cal D}';\mu_\lambda)$ such
that, for each $\xi\in{\cal D}$, the image of
$a_\lambda^\sim(\xi)$ under $I_\lambda$ is the operator of
multiplication by the function $\la\cdot,\xi\ra$ and
$I_\lambda\Omega=1$.

Let us dwell upon the explicit form of $I_\lambda$. Let
$\Ffin^*({\cal D})$ denote the dual space of $\Ffin({\cal D})$.
$\Ffin^*({\cal D})$ is the topological direct sum of the dual
spaces $({\cal D}_\C^{\hotimes n})^*$ of ${\cal D}_\C^{\hotimes
n}$. It will be convenient for us to realize each $({\cal
D}_\C^{\hotimes n})^*$ as the dual of ${\cal D}_\C^{\hotimes n}$
with respect to the zero space ${\cal H}_\C^{\hotimes n}\,n!$, so
that $({\cal D}_\C^{\hotimes n})^*$ becomes ${\cal
D}_\C^{\prime\,\hotimes n}$. Thus, $\Ffin^*({\cal D})$ consists of
infinite sequences $F=(F^{(n)})_{n=0}^\infty $, where $F^{(n)} \in
{\cal D }_\C^{\prime\,\hotimes n}$, and the dualization with
$f=(f^{(n)})_{n=0}^\infty\in\Ffin({\cal D})$ is given by
\begin{equation}\label{hxvghyf}\la\!\la F,f\ra\!\ra=\sum_{n=0}^\infty \la
\overline{F^{(n)}},f^{(n)} \ra\,n!,\end{equation} where
$\la\cdot,\cdot\ra$ denotes the dualization generated by the
scalar product in ${\cal H}^{\hotimes n}$, which is supposed to be
linear in  both dots.

 Next, according to the projection spectral
theorem, for $\mu_\lambda$-a.e.\ $\omega\in{\cal D}'$, there
exists a generalized joint vector
$P(\omega)=(P^{(n)}(\omega))_{n=0}^\infty\in \Ffin^*({\cal D})$ of
the family $a(\xi)$, $\xi\in{\cal D}$:
\begin{equation}\label{uhegfu} \forall f\in \Ffin({\cal D}):\qquad
\la\!\la P(\omega), a(\xi)f \ra\!\ra=\la \omega,\xi\ra\, \la\!\la
P(\omega), f \ra\!\ra,\end{equation} and for each
$f=(f^{(n)})_{n=0}^\infty\in \Ffin({\cal D})$ the action of
$I_\lambda$ onto $f$ is  given by
\begin{equation}\label{za67ftg} I_\lambda f=(I_\lambda
f)(\omega)= \la\!\la P(\omega), f \ra\!\ra =\sum_{n=0}^\infty \la
P^{(n)}(\omega),f^{(n)} \ra\,n!.\end{equation} We denote
$\wick{n}{:=}P^{(n)}(\omega)\,n!$, which is an element of ${\cal D
}_\C^{\prime\,\hotimes n}$ for $\mu_\lambda$-a.e.\ $\omega\in{\cal
D }'$. By \eqref{uhegfu} and \eqref{za67ftg}, we have, for
$\mu_\lambda$-a.e.\ $\omega\in{\cal D }'$,
\begin{gather} \la \omega,\xi\ra \la\wick{n},\xi^{\otimes
n}\ra=\la \wick{(n+1)},\xi^{\otimes(n+1)}\ra\notag\\\text{}+ \la
\wick{n},\lambda n(\xi^2)\hotimes
\xi^{\otimes(n-1)}+c_\lambda\la\xi\ra \xi^{\hotimes n}\ra\notag\\
\text{}+ \la \wick{(n-1)},
n\la\xi^2\ra\xi^{\otimes(n-1)}+n(n-1)(\xi^3)\hotimes
\xi^{\otimes(n-2)}\ra \label{zuaF}\end{gather} for all $\xi\in
{\cal D }$. Therefore, for $\mu_\lambda$-a.e.\ $\omega\in{\cal
D}'$, the $\wick{n}$'s are given by the recurrence relation
\eqref{udsH}. As easily seen, $\wick{n}$ is even well-defined as
an element of ${\cal D}_\C^{\prime\,\hotimes n}$ for each
$\omega\in{\cal D}'$.

Let us calculate the Fourier transform of $\mu_\lambda$. Let
$\Delta\in{\cal O}_c(X)$ and let $\chi_\Delta$ denote the
indicator of $\Delta$. One easily checks that each of the vectors
$\chi_\Delta^{\otimes n}$, $n\in\Z_+$, $\chi_\Delta^{\otimes
0}=\Omega$, belongs to $\Fext({\cal H})$, and let ${\cal
K}_\Delta$ be the subspace of $\Fext({\cal H})$ spanned by these
vectors.
 Let $a_\lambda^\sim(\chi_\Delta)$ denote the
operator in $\Fext({\cal H})$ whose image under the unitary
$I_\lambda$ is the operator of multiplication by the function
$$\la\omega,\chi_\Delta\ra{:=}\la \wick
1,\chi_\Delta\ra+c_\lambda\sigma(\Delta)=(I_\lambda\chi_\Delta)(\omega)+
c_\lambda\sigma(\Delta). $$ We approximate the indicator function
$\chi_\Delta$ by functions $\varphi_n\in{\cal D}$, $n\in\N$, such
that\linebreak  $\cup_{n\in\N}\supp\varphi_n$ is precompact in
$X$, $\varphi_n$'s are uniformly bounded, and $\varphi_n(x)\to
\chi_\Delta(x)$ as $n\to\infty$ for each $x\in X $. By using the
definitions of $\Fext({\cal H})$ and  $a^\sim_\lambda (\varphi)$
and the majorized convergence theorem, we get
$$a_\lambda^\sim(\chi_\Delta)\chi_\Delta^{\otimes n
}=\chi_\Delta^{\otimes(n+1)}+(\lambda n
+c_\lambda\sigma(\Delta))\chi_\Delta^{\otimes
n}+n(n-1+\sigma(\Delta))\chi_\Delta^{\otimes(n-1)}.$$ Therefore,
${\cal K}_\Delta$ is an invariant subspace for the operator
$a^\sim(\chi_\Delta)$. Hence, analogously to
\cite[pp.~315--316]{KL}, we may conclude that the distribution of
the random variable $\la\cdot,\chi_\Delta\ra$ is the probability
measure $\mu_{\lambda,\Delta}$ on $(\R,{\cal B}(\R))$ whose
orthogonal polynomials
$(P^{(n)}_{\lambda,\Delta}(x))_{n=0}^\infty$ with  leading
coefficient 1 satisfy the recurrence relation
\begin{equation}\label{iudyfv} x P^{(n)}_{\lambda,\Delta}(x)=
P^{(n+1)}_{\lambda,\Delta}(x)+(\lambda
n+c_\lambda\sigma(\Delta))P^{(n)}_{\lambda,\Delta}(x)
+n(n-1+\sigma(\Delta))P^{(n-1)}_{\lambda,\Delta}(x)\end{equation}
(the measure $\mu_{\lambda,\Delta}$ being defined uniquely through
this condition). Thus,
$(P^{(n)}_{\lambda,\Delta}(\cdot))_{n=0}^\infty$ is a system of
Laguerre polynomials for $\lambda=2$, and a  system of Meixner
polynomials for $\lambda\ne 2$. We will now consider only the case
$\lambda\ne2$, the case $\lambda=2$ being considered analogously
(see also \cite{KL}). By \cite{Meixner}, the Fourier transform of
the measure $\mu_{\lambda,\Delta}$ in a neighborhood of zero in
$\R$ is given by
\begin{equation}\label{gtzdtrd}\int_\R
e^{iux}\,\mu_{\lambda,\Delta}(dx)=\bigg(\frac{\alpha-\beta}{\alpha
e ^{-i \beta u }-\beta e^{-i\alpha
u}}\bigg)^{\sigma(\Delta)/(\alpha\beta)}\exp[ic_\lambda
u\sigma(\Delta)].\end{equation} Therefore, for any $u\in\R$
satisfying \begin{equation}\label{sjihgu} \bigg| \frac{\alpha
(e^{-i\beta u }-1)-\beta(e^{i\alpha u}-1)
}{\alpha-\beta}\bigg|<1,\end{equation} we get
\begin{gather} \int_{{\cal D}'}
e^{iu\la\omega,\chi_{\Delta}\ra}\,\mu_\lambda(d\omega)=\exp\bigg[
-\frac{\sigma(\Delta)}{\alpha\beta}\log\bigg(\frac{\alpha
e^{-i\beta u}-\beta e^{-i\alpha u}}{\alpha-\beta}\bigg)
+ic_\lambda u\sigma(\Delta)\bigg]\notag\\ =\exp\bigg[
-\frac1{\alpha\beta}\int_X\log\bigg(\frac{\alpha e^{-i\beta u
\chi_\Delta(x)}-\beta e^{-i\alpha
u\chi_\Delta(x)}}{\alpha-\beta}\bigg)\,\sigma(dx)+i
c_\lambda\int_X u\chi_\Delta(x) \,\sigma(dx)\bigg].\label{hiadfcu}
\end{gather}

Now, let $\Delta_1,\dots,\Delta_n\in{\cal O}_c(X)$ be disjoint.
Then the spaces ${\cal K}_{\Delta_1}\ominus {\cal F}_0({\cal
H})$,\dots, ${\cal K}_{\Delta_n}\ominus {\cal F}_0({\cal H})$ are
orthogonal in $\Fext({\cal H})$. Therefore, the random variables
$\la\cdot,\chi_{\Delta_1}\ra,\dots,\la\cdot,\chi_{\Delta_n}\ra$
are independent. Hence, for any step function
$\varphi=\sum_{i=1}^n u_i\chi_{\Delta_i}$ such that all $u_i$'s
satisfy \eqref{sjihgu} with $u=u_i$, formula \eqref{hisav} holds.

Finally, fix any $\varphi\in{\cal D}$ satisfying \eqref{hisav}.
Choose any sequence of step functions $\{\varphi_n\}_{n\in\N}$
such that $$ \sup_{n\in\N,\, x\in X}\bigg| \frac{\alpha
(e^{-i\beta\varphi_n(x)}-1)-\beta(e^{i\alpha\varphi_n(x)}-1)
}{\alpha-\beta}\bigg|<1,$$ $\cup_{n\in\N}\operatorname{supp}
\varphi_n\in{\cal O}_c(X)$, and $\varphi_n$ converges pointwisely
to $\varphi$ as $n\to\infty$.  Then, by the majorized convergence
theorem, we conclude that the right hand side of \eqref{hisav}
with $\varphi=\varphi_n$ converges to the right hand side of
\eqref{hisav}. On the other hand, $\la\cdot,\varphi_n\ra$
converges to $\la\cdot,\varphi\ra$ in $L^2({\cal D}';\mu)$, and
therefore also in probability. Hence, again by the majorized
convergence theorem, the left hind side of \eqref{hisav} with
$\varphi=\varphi_n$ converges to the left hand side of
\eqref{hisav}.\quad $\blacksquare$

\begin{corollary}\label{fsopioipioi}For each $\Delta\in{\cal O}_c(X)$, the
distribution $\mu_{\lambda,\Delta}$ of the random variable
$\la\cdot,\chi_\Delta\ra$ under $\mu_\lambda$ is given as
follows\rom: For $\lambda>2$\rom, $\mu_{\lambda,\Delta}$ is the
negative binomial \rom(Pascal\rom) distribution
\begin{equation}\label{formula1}
\mu_{\lambda,\Delta}=(1-p_{\lambda})^{\sigma(\Delta)}\sum_{k=0}^\infty
 \frac{\big(\sigma(\Delta)\big)_k}{k!}\, p_\lambda^k \, \delta_{\sqrt{\lambda^2-4}\,k}, \end{equation} where $$
p_\lambda{:=}\frac{\lambda-\sqrt{\lambda^2-4}}{\lambda+\sqrt{\lambda^2-4}},$$
and for $r> 0$ $(r)_0{:=}1$, $(r)_k{:=}r(r+1)\dotsm (r+k-1)$\rom,
$k\in\N$\rom. For $\lambda=2$, $\mu_{2,\Delta}$ is the Gamma
distribution
\begin{equation}\label{formula2}\mu_{2,\Delta}(ds)=\frac{s^{\sigma(\Delta)-1}e^{-s}}{\Gamma(\sigma(\Delta))}\,\chi_{[0,\infty)}(s)\,
ds.\end{equation} Finally\rom, for $\lambda\in[0,2)$\rom,
\begin{multline}\label{formula3}\mu_{\lambda,\Delta}(ds)=\frac{(4-\lambda^2)^{(\sigma(\Delta)-1)/2}}{2\pi
\Gamma(\sigma(\Delta)) }\\ \times \big|
\Gamma\big(\sigma(\Delta)/2+i(4-\lambda^2)^{-1/2}s\big)
\big|^2\,\exp\big[
-s2(4-\lambda^2)^{-1/2}\arctan\big(\lambda(4-\lambda^2)^{-1/2}\big)
\big]\,ds.\end{multline}
\end{corollary}

\noindent {\it Proof}. The result follows from (the proof of)
Theorem~\ref{z8fzvear},  \cite{Meixner}, and \cite[Ch.~VI,
sect.~3]{Chihara}
 (see also \cite{ST} and \cite[subsec.~4.3.5]{S}).\quad $\blacksquare$

For $f^{(n)}\in\Fext^{(n)}({\cal H})$\rom, let
$\la{:}\cdot^{\otimes n}{:}_\lambda,f^{(n)}\ra$ denote the element
of $L^2({\cal D}';\mu_\lambda)$ defined as $I_\lambda f^{(n)}$.

\begin{corollary}\label{tzasdfafv} For any  $\Delta\in{\cal O}_c(X)$\, we have
\begin{equation}\label{ijhft}
\la \wick n,\chi_\Delta^{\otimes
n}\ra=P^{(n)}_{\lambda,\Delta}(\la\omega,\chi_\Delta\ra)\qquad
\text{$\mu_\lambda$-a.a.\ $\omega\in{\cal D}'$,}\end{equation}
where $(P_{\lambda,\Delta}^{(n)})_{n=0}^\infty$ are the
polynomials on $\R$ as in the proof of
Theorem~\rom{\ref{z8fzvear}}\rom.\end{corollary}

\noindent {\it Proof}. This result directly follows from the proof
of Theorem~\ref{z8fzvear}.\quad $\blacksquare$

\begin{remark}\label{bhgfzt}\rom{ Let us  state the one-dimensional analog
of the results of Theorem~\ref{z8fzvear} and
Corollaries~\ref{fsopioipioi}, \ref{tzasdfafv}. We consider the
weighted $\ell_2$-space ${\cal F}=\ell_2(((n!)^2)_{n=0}^\infty)$
consisting of all sequences $f=(f^{(n)})_{n=0}^\infty$,
$f^{(n)}\in\C$, such that $\|f\|_{\cal F}^2{:=}\sum_{n=0}^\infty
|f^{(n)}|^2 (n!)^2<\infty$. Let ${\cal F}_{\mathrm fin}$ denote
the set of all finite sequences in $\cal F$. For each
$\lambda\in[0,\infty)$, we define a linear Hermitian  operator
$a_\lambda$ in $\cal F$ with domain ${\cal F}_{\mathrm fin}$ by
setting $$ a_\lambda =a^+ +\lambda a^0 + a^-+ c_\lambda
\operatorname{id},$$ where
$(a^+f^{(n)})^{(k)}=\delta_{k,n+1}\,f^{(n)}$, $(a^0
f^{(n)})^{(k)}=\delta_{k,n}\, nf^{(n)}$,
$(a^-f^{(n)})^{(k)}=\delta_{k,n-1}\,n^2f^{(n)}$. Here
$\delta_{i,j}=1$ if $i=j$ and $\delta_{i,j}=0$ otherwise, and
$c_\lambda$ is given by \eqref{qwqwuzzuuz}. We note that, under
the natural unitary mapping of the weighted $\ell_2$-space $\cal
F$ onto the usual $\ell_2$, the operator $a_\lambda$ goes over
into the operator defined by the infinite Jacobi matrix
$J=(\alpha_{m,n})_{m,n=0}^\infty$ with the elements
$\alpha_{n,n}=\lambda n+c_\lambda$, $n\in\Z_+$,
$\alpha_{n,n+1}=\alpha_{n+1,n}=n+1$, $n\in\Z_+$, and
$\alpha_{m,n}=0$ for $|m-n|>1$.

The operator $a_\lambda$ is essentially self-adjoint, and let
$a_\lambda^\sim$ denote its closure. By  the spectral theory of
infinite Jacobi matrices (e.g.\ \cite[Ch.~VII, Sect.~1]{b}), there
exist a unique probability measure ${\frak m}_\lambda$ on
$(\R,{\cal B}(\R))$ and a unique unitary operator ${\cal
I}_\lambda:{\cal F}\to L^2(\R;{\frak m}_\lambda)$ such that the
image of the operator $a_\lambda^\sim$ under ${\cal I}_\lambda$ is
the operator of multiplication by the variable $x$ and ${\cal
I}_\lambda(1,0,0,\dots)=1$. The mapping ${\cal I}_\lambda$ is
given on the dense set ${\cal F}_{\mathrm fin}$ by $${\cal
F}_{\mathrm fin}\ni f=(f^{(n)})_{n=0}^\infty\mapsto {\cal
I}_\lambda f=\sum_{n=0}^\infty f^{(n)}P_\lambda^{(n)}(x).$$ Here,
$(P^{(n)}_\lambda)_{n=0}^\infty$ is the system of polynomials on
$\R$ satisfying the following recurrence relation: for $n\in\Z_+$
$$ xP_\lambda^{(n)}(x)=P_\lambda^{(n+1)}(x)+(\lambda
n+c_\lambda)P_\lambda^{(n)}(x)+ n^2 P_\lambda^{(n-1)}(x),\qquad
P_\lambda^{(-1)}(x)=0,\ P_\lambda^{(1)}(x)=1,$$ and ${\frak
m}_\lambda$ is the unique probability measure on $\R$ with respect
to which the polynomials $(P^{(n)}_\lambda)_{n=0}^\infty$ are
orthogonal. By \cite{Meixner}, the Fourier transform of the
measure ${\frak m }_\lambda$ in a neighborhood of zero in $\R$ is
given by
\begin{equation}\label{tzrdt} \int_\R
e^{iux}\,{\frak m}_\lambda(dx)=\begin{cases} e^{ic_\lambda
u}\left(\frac{\alpha-\beta}{\alpha e^{-i\beta u}-\beta e^{-i\alpha
u}}\right)^{1/(\alpha\beta)},& \lambda\ne 2,\\
(1-iu)^{-1},&\lambda=2.\end{cases}\end{equation} Furthermore, the
measure ${\frak m}_\lambda$ is explicitly given by the right hand
side of
 formula \eqref{formula1}, resp.\ \eqref{formula2}, resp.\ \eqref{formula3}, with
$\sigma(\Delta)=1$. Thus ${\frak m}_\lambda$ is a Pascal
distribution for $\lambda\in[0,2)$, a Gamma distribution for
$\lambda=2$, and a Meixner distribution for $\lambda>2$.

}\end{remark}

 We will now show that each $\mu_\lambda$ is a compound Poisson, respectively L\'evy noise
 measure.

\begin{corollary} \label{audstz} For each $\lambda\ge2$,
$\mu_\lambda$ is  a compound Poisson measure on $({\cal D}',{\cal
C }_\sigma({\cal D}'))$ whose  L\'evy--Khintchine representation
of the Fourier transform reads as follows\rom:
\begin{equation}\label{hjsarfguzwegr}
\int_{{\cal D
}'}e^{i\la\omega,\varphi\ra}\,\mu_\lambda(d\omega)=\exp\bigg[\int_{X\times\R_+}(e^{is\varphi(x)}
-1)\,\sigma(dx)\,\nu_\lambda(ds)\bigg],\qquad \varphi\in{\cal
D},\end{equation} where \begin{align}
\nu_2(ds)&=\frac{e^{-s}}{s}\,ds,\notag\\
\nu_\lambda(ds)&=\sum_{k=1}^\infty
\frac{p_\lambda^k}{k}\,\delta_{\sqrt{\lambda^2-4}\,k},\qquad
\lambda>2.\label{jkwarf}\end{align} In particular\rom, each
$\mu_\lambda$ is concentrated on the set of all Radon measures on
$(X,{\cal B }(X))$ of the form $\sum_{n=1}^\infty s_n
\delta_{x_n}$, $\{x_n\}\subset X$\rom, $s_n>0$\rom, $n\in\N$\rom.
\end{corollary}

\noindent{\it Proof}. It follows from the general theory of
compound Poisson measures (e.g.\ \cite{Kal83}) that there exists a
compound Poisson measure $\tilde\mu_\lambda$ whose Fourier
transform is given by \eqref{hjsarfguzwegr} with  $\nu_\lambda$
given by \eqref{jkwarf}, and which is concentrated on the set of
those Borel measures on $X$ as in the formulation of the theorem.
Furthermore, it follows from the general theory that, for any
disjoint $\Delta_1,\dots,\Delta_n\in{\cal O}_c(X)$, the random
variables are independent. Thus, it suffices to show that, for any
fixed $\Delta\in {\cal O}_c(X)$, the distributions of the random
variable $\la\omega,\chi_\Delta\ra$ under $\mu_\lambda$ and
$\tilde \mu_\lambda$ coincide. But this can be easily done by
calculating the Fourier transform $$ \int_{{\cal D}'}e^{i u
\la\omega,\chi_\Delta\ra}\,\tilde\mu_\lambda(d\omega)=
\exp\bigg[\sigma(\Delta)\int_{\R_+}(e^{ius}-1)\,\nu_\lambda(ds)
\bigg],\qquad u\in\R,$$ and comparing it with the Fourier
transform of the measure $\mu_{\lambda,\Delta}$ (see
\eqref{gtzdtrd} for the case $\lambda>2$).\quad$\blacksquare$

In the case $\lambda\in[0,2)$, the situation is a little  bit more
complicated, since the corresponding L\'evy  measure $\nu_\lambda$
to be identified  does not have the first moment finite.

\begin{corollary}\label{uiawfruza}
For $\lambda\in[0,2)$\rom, $\mu_\lambda$ is the L\'evy  noise
measure whose  L\'evy--Khintchine representation of the Fourier
transform reads as follows\rom:
\begin{equation}\label{juihscu} \int_{{\cal D}'}
e^{i\la\omega,\varphi\ra}\,\mu_\lambda(d\omega)=\exp\bigg[\int_{X\times\R}\big(e^{is\varphi(x)}
-1-is\varphi(x)\big)\,\sigma(dx)\,\nu_\lambda(ds)+ic_\lambda
\la\varphi\ra \bigg],\qquad \varphi\in{\cal D},\end{equation}
where
\begin{multline}\label{iowgzs}
\nu_\lambda(ds)= \frac{\sqrt{4-\lambda^2}}{2\pi}\\ \times
\big|\Gamma\big(1+i
(4-\lambda^2)^{-1/2}s\big)\big|^2\,\exp\big[-s2(4-\lambda^2)^{-1/2}\arctan
\big(\lambda(4-\lambda^2)^{-1/2}\big) \big]\,\frac1{s^2} \,ds.
\end{multline}

\end{corollary}

\noindent{\it Proof}. The existence of a probability measure
$\tilde\mu_\lambda$ on ${\cal D}'$ whose Fourier transform is
given by the right hand side of \eqref{juihscu} with $\nu_\lambda
$ given by \eqref{iowgzs} follows by e.g.\ the Bochner--Minlos
theorem. Furthermore, as easily seen, for each $\Delta\in{\cal
O}_c(X)$, one can naturally define a random variable
$\la\cdot,\chi_\Delta\ra$ as an element of $L^2({\cal D
}';\tilde\mu_\lambda)$. By \eqref{juihscu}, for any disjoint
$\Delta_1,\dots,\Delta_n\in{\cal O}_c(X)$, the random variables
$\la\cdot,\chi_{\Delta_1}\ra,\dots,\la\cdot,\chi_{\Delta_n}\ra$
are independent. Analogously to the proof of
Corollary~\ref{audstz}, we conclude the statement by calculating
the integral $\int_\R(e^{ius}-1-ius)\,\nu_{\lambda}(ds)$,
$u\in\R$, using
 \cite[subsec.~4.3.5]{S}, see also \cite{ST}.\quad$\blacksquare$

\begin{remark}\rom{It follows from Corollaries~\ref{audstz} and \ref{uiawfruza}
 that $s^2\nu_\lambda(ds)$ is a Meixner distribution for
$\lambda\in[0,2)$, gamma distribution for $\lambda=2$, and Pascal
distribution for $\lambda>2$. Furthermore, for each $\lambda\ge0$,
$s^2\,\nu_\lambda(ds)$ is a probability measure on $\R$ whose
orthogonal polynomials $(Q^{(n)}_{\lambda})_{n=0}^\infty$ with
leading coefficient 1 satisfy the following recurrence relation:
$$ sQ^{(n)}_{\lambda}(s)= Q^{(n+1)}_\lambda(s)+ \lambda
(n+1)Q^{(n)}(s)+n(n+1)Q^{(n-1)}_\lambda(s),\qquad n\in \Z_+,\
Q^{-1}_\lambda(s){:=}0. $$
 }\end{remark}

We denote by ${\cal P}({\cal D}')$ the set of continuous
polynomials on ${\cal D}'$, i.e., functions on ${\cal D}'$ of the
form $$F(\omega)=\sum_{i=0}^n\la\omega^{\otimes
i},f^{(i)}\ra,\qquad f^{(i)}\in{\cal D}_\C^{\hotimes i},\
\omega^{\hotimes 0}{:=}1,\ i\in\Z_+. $$ The greatest number $i$
for which $f^{(i)}\ne0$ is called the power of a polynomial.  We
denote by ${\cal P}_n({\cal D}')$ the set of continuous
polynomials of power $\le n$.

\begin{corollary}
\label{uiwef} For each $\lambda\ge 0$\rom, we have
$I_\lambda(\Ffin({\cal D}))={\cal P}({\cal D }')$\rom. In
particular\rom, ${\cal P}({\cal D}')$ is a dense subset of
$L^2({\cal D}';\mu_\lambda)$\rom. Furthermore\rom, let ${\cal
P}_{\lambda,n}({\cal D}')$ denote the closure of ${\cal P}_n({\cal
D}')$ in $L^2({\cal D}';\mu_\lambda)$\rom, and let
$(L^2_{\lambda,n})$ denote the orthogonal difference ${\cal
P}_{\lambda,n}({\cal D}')\ominus {\cal P}_{\lambda,n-1}({\cal
D}')$ in L$^2({\cal D}';\mu_\lambda)$\rom. Then\rom,
\begin{equation}\label{kjft} L^2({\cal D}';\mu_\lambda)=\bigoplus_{n=0}^\infty
(L^2_{\lambda,n}).\end{equation} Finally, let $P_{\lambda,n}$
denote the orthogonal  projection of $L^2({\cal D}';\mu_\lambda)$
onto $(L^2_{\lambda,n})$\rom. Then\rom, for any $f^{(n)}\in{\cal
D}_\C^{\hotimes n}$\rom, \begin{equation}\label{hzdtrds}
P_{\lambda,n}(\la\cdot^{\otimes
n},f^{(n)}\ra)=\la{:}\cdot^{\otimes n}{:}_\lambda,f^{(n)}\ra\qquad
\text{$\mu_\lambda$-a\rom.e\rom.} \end{equation} and
\begin{equation}\label{jhafd} I_\lambda(\Fext^{(n)}({\cal
H}))=(L^2_{\lambda, n}).\end{equation}
\end{corollary}

\noindent{\it Proof}. Using recurrence  relation \eqref{udsH}, we
obtain by induction the inclusion  $I_\lambda(\Ffin({\cal
D}))\subset{\cal P}({\cal D}')$ and moreover, the equality
\begin{equation}\label{hjsvag}
\la\wick{n},f^{(n)}\ra=\la\omega^{\otimes
n},f^{(n)}\ra+p_{n-1}(\omega),\qquad f^{(n)}\in{\cal
D}_\C^{\hotimes n},
\end{equation}
where  $p_{n-1}\in{\cal P}_{n-1}({\cal D}')$. Using
\eqref{hjsvag}, we then  obtain by induction also the inverse
inclusion ${\cal P}({\cal D}')\subset I_\lambda( \Ffin({\cal
D}'))$. That ${\cal P}({\cal D}')$ is dense in $L^2({\cal
D}';\mu_\lambda)$ follows from the fact that $\Ffin( {\cal D})$
 is dense in $\Fext({\cal H})$. Decomposition  \eqref{kjft} now
 becomes evident. Finally, \eqref{hzdtrds} follows by
 \eqref{hjsvag}, and \eqref{jhafd} is a consequence
 of \eqref{hzdtrds}.\quad $\blacksquare$

\section{The generating function}\label{jfrtd}

Now, we will identify the generating function of the polynomials
$\la\wick n,\varphi^{\otimes n}\ra$, $\varphi\in{\cal D}$.

 Let
${\cal M}(X)$ denote the set of signed Radon measures on $(X,{\cal
B }(X))$. We can evidently identify any measure $\omega\in{\cal
M}(X)$ with an element $\tilde\omega\in{\cal D}'$ by setting
$$\la\tilde\omega,\varphi\ra{:=}\int _X
\varphi(x)\,\omega(dx),\qquad \varphi\in{\cal D}.$$ In what
follows, we will just  write $\omega$ instead of $\tilde\omega$.
Then, for $\Delta\in{\cal O }_c(X)$, the writing
$\la\omega,\chi_\Delta\ra$ will  mean $\omega(\Delta)$.

\begin{proposition}\label{ihaugcuzd}We have\rom, for $\lambda\ne 2$\rom,
\label{uizt}
\begin{multline}\label{hjdtvgr}
G_\lambda(\omega,\varphi){:=}\sum_{n=0}^\infty
\frac1{n!}\,\la\wick n,\varphi^{\otimes
n}\ra=\exp\bigg[-\frac1{\alpha-\beta}
\bigg\la\log\bigg(\frac{(1-\beta\varphi)^{1/\beta}}{(1-\alpha\varphi)^{1/\alpha}}\bigg)\bigg\ra
\\ \text{}+\frac1{\alpha-\beta}\bigg\la\omega-c_\lambda,\log\bigg(\frac{1-\beta\varphi}{1-\alpha\varphi}\bigg)\bigg\ra
\bigg]
\end{multline} and for $\lambda=2$
\begin{equation}\label{ioutz}G_2(\omega,\varphi)
{:=}\sum_{n=0}^\infty \frac1{n!}\,\la {:}\omega^{\otimes
n}{:}_2,\varphi^{\otimes n}\ra=\exp\bigg[-\la
\log(1+\varphi)\ra+\bigg\la\omega,\frac\varphi{\varphi+1}\bigg\ra\bigg].\end{equation}
Formulas \eqref{hjdtvgr}\rom, \eqref{ioutz} hold for each
$\omega\in{\cal M}(X)$ and for each $\varphi\in{\cal D}_\C$
satisfying $\|\varphi\|_\infty<(\max(|\alpha|,|\beta|))^{-1}$ for
\eqref{hjdtvgr} and $\|\varphi\|_\infty<1$ for \eqref{ioutz}\rom.
More generally\rom, for each fixed $\tau\in T$\rom, there exists a
neighborhood  of zero in ${\cal D}_\C$ \rom(depending on
$\lambda$\rom)\rom, denoted by ${\cal O}_\tau$\rom, such that
\eqref{hjdtr}\rom, respectively \eqref{ioutz}, holds for all
$\omega\in{\cal H}_{-\tau}$ and all $\varphi\in  {\cal
O}_\tau$\rom.
\end{proposition}

\begin{remark}\label{fzfrzdzrz}\rom{In the one-dimensional case (see
Remark~\ref{bhgfzt}), the generating function of the polynomials
$(P^{(n)}_\lambda(\cdot))_{n=0}^\infty$ is given by (cf.\
\cite{Meixner})
\begin{align*}{\cal G}_\lambda(x,u){:=}\sum_{n=0}^\infty \frac{u^n}{n!}\,P^{(n)}_\lambda (x)&=\left(
\frac{(1-\beta u)^{1/\beta}}{(1-\alpha
u)^{1/\alpha}}\right)^{-1/(\alpha-\beta)} \left( \frac{1-\beta
u}{1-\alpha u}\right)^{(x-c_\lambda)/(\alpha-\beta)},\qquad
\lambda\ne 2,\\ {\cal G}_2(x,u){:=}\sum_{n=0}^\infty
\frac{u^n}{n!}\,P^{(n)}_2 (x)&= \frac{1}{1+u}\,
e^{ux/(u+1)}\end{align*} for $u$ from a neighborhood of zero in
$\C$. }\end{remark}

\noindent{\it Proof of Proposition~\rom{\ref{ihaugcuzd}}}. We only
prove formula \eqref{hjdtvgr}, corresponding to the case
$\lambda\ne 2$. Let us fix $\omega\in{\cal M}(X)$. Then, as easily
seen from \eqref{udsH}, $\wick n\in{\cal M}(X^n)$ for each
$n\in\N$ (if we identify ${\cal M}(X^n)$ as a subset of ${\cal
D}^{\prime\,\otimes n})$. Fix $\Delta\in{\cal O}_c(X)$. As follows
from the proof of Theorem~\ref{z8fzvear}, the equality in
\eqref{ijhft} holds for each $\omega\in{\cal M}(X)$.  Then,
$$\sum_{n=0}^\infty\frac{1}{n!}\,\la\wick{n},(u\chi_{\Delta})^{\otimes
n }\ra=
\sum_{n=0}^\infty\frac{u^n}{n!}\,\la\wick{n},\chi_{\Delta}^{\otimes
n }\ra=\sum_{n=0}^\infty\frac{u^n}{n!}\,
P_{\lambda,\Delta}^{(n)}(\la\omega,\chi_{\Delta}\ra). $$  Hence,
it follows from  \cite{Meixner} (see also Remark~\ref{fzfrzdzrz})
that \eqref{hjdtr} holds with $\varphi=u\chi_\Delta$ and $u\in\C$
such that $|u|<(\max(|\alpha|,|\beta|))^{-1}$.

We next prove the following lemma.

\begin{lemma}\label{sjihd} For any $\omega\in{\cal M}(X)$ and any
disjoint $\Delta_1,\dots,\Delta_l\in{\cal O}_c(X)$\rom,
\begin{equation}\label{zftd}\la\wick n,\chi_{\Delta_1}^{\otimes
k_1}\hotimes\dotsm\hotimes \chi_{\Delta_l}^{\otimes
k_l}\ra=\prod_{i=1}^l\la \wick{k_i},\chi_{\Delta_i}^{\otimes
k_i}\ra =\prod_{i=1}^l
P_{\lambda,\Delta}^{(k_i)}(\la\omega,\chi_{\Delta_i}\ra),\end{equation}
where $k_1,\dots,k_l\in\N$\rom, $k_1+\dots+k_l=n$\rom.
\end{lemma}

\noindent{\it Proof}. We prove this lemma by induction in
$n\in\N$. For $n=1$, formula \eqref{zftd} trivially holds. Now,
suppose that \eqref{zftd} holds for all $n\le N$. Let
$k_1,\dots,k_l\in\N$, $k_1+\dots+k_l=N+1$. Applying  recurrence
formula \eqref{udsH} to $\wick{(N+1)}$, we express $\la
\wick{(N+1)},\chi_{\Delta_1}^{\otimes k_1}\hotimes\dotsm\hotimes
\chi_{\Delta_l}^{\otimes k_l}\ra$ through
$\la\omega,\chi_{\Delta_j}\ra$, $\la \wick
N,\chi_{\Delta_1}^{\otimes k_1}\hotimes\dotsm\hotimes
\chi_{\Delta_j}^{\otimes (k_j-1)}\hotimes\dotsm\hotimes
\chi_{\Delta_l}^{\otimes k_l}\ra$, $\la \wick
{(N-1)},\chi_{\Delta_1}^{\otimes k_1}\hotimes\dotsm\hotimes
\chi_{\Delta_j}^{\otimes (k_j-2)}\hotimes\dotsm\hotimes
\chi_{\Delta_l}^{\otimes k_l}\ra$, and $\sigma(\Delta_j)$,
$j=1,\dots,l$. Applying  formula \eqref{zftd} with $n=N$ and
$n=N-1$ and then using the recurrence relation \eqref{iudyfv} for
the polynomials $P^{(n)}_{\lambda,\Delta}$, we conclude the
statement.\quad $\blacksquare$

Fix any disjoint $\Delta_1,\dots,\Delta_m\in{\cal O}_c(X)$ and any
$u_1,\dots,u_m\in\C$ satisfying
$|u_i|<\max(|\alpha|,|\beta|)^{-1}$, $i=1,\dots,m$, and set
$f{:=}\sum_{i=1}^m u_i\chi_{\Delta_i}$. By Lemma~\ref{sjihd}, we
get $$\sum_{n=0}^\infty\frac1{n!}\, \la\wick n,f^{\otimes
n}\ra=\prod_{i=1}^m\bigg( \sum_{n=0}^\infty\frac{u_i^n}{n!}\,\la
\wick n,\chi_{\Delta_i}^ {\otimes n}\ra\bigg).$$ Hence,
\eqref{hjdtr} holds with $\varphi=f$.

Using  \eqref{udsH}, one easily shows by induction that, for each
fixed  $\Lambda\in{\cal O}_c(X)$, there exists a constant
$C_{\Lambda,\omega}$ such that \begin{equation}\label{ghxd}
\forall n\in\N:\qquad \big|\wick n \restriction \Lambda^n\big| \le
n!\, C_{\Lambda,\omega}^n,
\end{equation}
where $\big|\wick n \restriction \Lambda^n\big|$ denotes the full
variation of the signed measure $\wick n$ on $\Lambda^n$. Fix any
$\varphi\in{\cal D}$ such that  $\operatorname{supp}
\varphi\subset\Lambda$ and $\|\varphi\|_\infty<
C_{\Lambda,\omega}^{-1}$. Let $\{f_k,\,k\in\N\}$ be a sequence of
step functions on $X$ such that $$C{:=}\sup_{k\in\N,\, x\in
X}|f_k(x)|< C_{\Lambda,\omega}^{-1},$$
$\cup_{k\in\N}\operatorname{supp}f_k\subset\Lambda$ and $f_k\to
\varphi$ as $k\to\infty$ uniformly on $X$. We then get by
\eqref{ghxd}  \begin{align}& \bigg| \sum_{n=0}^\infty \frac1{n!}\,
\la\wick n, f_k^{\otimes n}\ra-\sum_{n=0}^\infty \frac{1} {n!}\,
\la \wick n,\varphi^{\otimes n}\ra\bigg|\notag\\ &\qquad \le
\sum_{n=0}^\infty \frac1{ n!}\,\big|\wick n \restriction
\Lambda^n\big|\sup_{(x_1,\dots,x_n)\in X^n}|f^{\otimes
n}_k(x_1,\dots,x_n)-\varphi^{\otimes n}(x_1,\dots,x_n)|\notag\\
&\qquad\le  \sum_{n=0}^\infty C_{\Lambda,\omega}^n \,n
\max(\|\varphi\|_\infty,C)^{n-1}\sup_{x\in
X}|f_k(x)-\varphi(x)|\to0\quad\text{as }k\to\infty.\label{hay}
\end{align}
Let $\widetilde G_\lambda(\omega,\varphi)$ denote the right hand
side of \eqref{hjdtr}. Then, if
$$\max(\|\varphi\|_\infty,C)<\max(|\alpha|,|\beta|)^{-1},$$ by the
majorized convergence theorem, $\widetilde G_\lambda(\omega,
f_k)\to \widetilde G_\lambda(\omega,\varphi)$ as $k\to\infty$.
Thus, \eqref{hjdtr} holds for any $\varphi\in{\cal D}$ such that
$\operatorname{supp}\varphi\subset \Lambda$ and
$$\|\varphi\|_\infty<\min(C_{\Lambda,\omega}^{-1},\max(|\alpha|,|\beta|)^{-1}).$$

Let us show that \eqref{hjdtr} holds  for any $\varphi\in{\cal D}$
such that $\|\varphi\|_\infty<\max(|\alpha|,|\beta|)^{-1}$. Fix
such $\varphi\in{\cal D}$. Denote $$
O_\varphi{:=}\big\{z\in\C:|z|<
1+(1/2)(\max(|\alpha|,|\beta|)^{-1}-\|\varphi\|_\infty) \big\},$$
and consider the analytic function $$ O_\varphi \ni z\mapsto
g(z){:=}\widetilde G_\lambda(\omega,z\varphi).$$ For all $z\in\C$
such that $|z|\le1$, we have the Taylor expansion
\begin{equation}\label{hfdtr} g(z)=\sum_{n=0}^\infty
\frac{g^{(n)}(0)}{n!}\, z^n.\end{equation} But it follows from the
proved above that \begin{equation}\label{uzdrt}
g(z)=\sum_{n=0}^\infty \frac{1}{n!}\,\la\wick n,\varphi^{\otimes
n}\ra z^n\end{equation} for all $z\in\C$ from some neighborhood of
zero. Comparing \eqref{hfdtr} and \eqref{uzdrt}, we get
\begin{equation}\label{gzzddtr} g^{(n)}(0)=\la\wick n,\varphi^{\otimes n}\ra,\end{equation} and
thus, by \eqref{hfdtr} and \eqref{gzzddtr}  we have proved the
proposition in the case where $\omega\in{\cal M}(X)$.

Let us consider the general case. Fix $\tau\in T$. Without loss of
generality, we can suppose that the inclusion ${\cal
H}_\tau\hookrightarrow L^2(X;\sigma)$ is quasi-nuclear,
$$\|\varphi\|_{\infty}\le C_\tau\|\varphi\|_\tau,\qquad
\varphi\in{\cal D},\ C_\tau>0, $$ and $1\in{\cal H}_{-\tau}$. Let
$$ {\cal O}_{\tau}{:=}\big\{\varphi\in {\cal D}_\C:
\|\varphi\|_\tau\le(2C_\tau \max(|\alpha|,|\beta|))^{-1} \big\}.$$
It follows from \eqref{hjdrt} that $$\sup_{\varphi\in{\cal
O}_\tau}\max\big(\|\log(1-\alpha\varphi)\|_\tau,\|\log(1-\beta\varphi)\|_\tau\big)<\infty.$$
Then, for each fixed $\omega\in{\cal D}_{-\tau}$, the function
$\widetilde G_\lambda(\omega,\cdot)$ is G-holomorphic and bounded
on ${\cal O}_\tau$. Hence, by e.g.\ \cite{Dineen}, $\widetilde
G_\lambda(\omega,\cdot)$ is holomorphic on ${\cal O}_\tau$. The
Taylor decomposition of $\widetilde G_\lambda(\omega,\cdot)$ and
the kernel theorem (cf.\ \cite[subsec.~4.1]{KSWY}) yield
\begin{equation}\label{ugfdrt} \widetilde G_\lambda(\omega,\varphi)=\sum_{n=0}^\infty \frac{1}{n!}\, \la
G_\lambda^{(n)}(\omega),\varphi^{\otimes n}\ra, \end{equation}
where $G_\lambda^{(n)}(\omega)\in{\cal H}_{-\tau}^{\hotimes n}$ is
given through $$\la G_\lambda^{(n)}(\omega),\varphi^{\otimes n}\ra
=\frac{d^n}{dt^n}\Big|_{t=0} \widetilde
G_\lambda(\omega,t\varphi).$$ Next, for each $\omega\in{\cal
M}(X)\cap {\cal H}_{-\tau}$, we have, by the proved above,
\begin{equation}\label{zadfoiv}\la
G_\lambda^{(n)}(\omega),\varphi^{\otimes n}\ra=\la\wick
n,\varphi^{\otimes n}\ra,\qquad \varphi\in{\cal D}.\end{equation}
Differentiating $\widetilde G_\lambda(\omega,t\varphi)$ in $t$, we
conclude that $\frac{d^n}{dt^n}\big|_{t=0}\widetilde
G_\lambda(\omega,t\varphi)$ depends continuously on
$\omega\in{\cal H_{-\tau}}$. On the other hand, $\la\wick
n,\varphi^{\otimes n}\ra$ does also depend continuously on
$\omega\in{\cal H_{-\tau}}$. Since ${\cal M}(X)\cap {\cal
H}_{-\tau}$ is dense in ${\cal H}_{-\tau}$ (${\cal M}(X)\cap {\cal
H}_{-\tau}$, in particular, contains $\cal D$), we conclude that
\eqref{zadfoiv} holds for all $\omega\in{\cal H}_{-\tau}$, and
hence also   \eqref{ugfdrt} holds for all $\omega\in{\cal
H}_{-\tau}$.\quad $\blacksquare$

\begin{corollary}\label{hxxd}For each $\lambda\ge0$\rom, the
function $$ \varphi\mapsto L_\lambda(\varphi){:=} \int_{{\cal D}'}
e^{\la\omega,\varphi\ra}\,\mu_\lambda(d\omega) $$ is well defined
and holomorphic on some neighborhood of zero in ${\cal D}_\C$.
Furthermore\rom, for each fixed $\tau\in T$, there exists a
neighborhood of zero in ${\cal D}_\C$, denoted by $ O_\tau$\rom,
such that\rom, for all $\omega\in{\cal H}_{-\tau}$ and all
$\varphi\in O_\tau$\rom, we have
\begin{equation}\label{ihugayzc}
G_\lambda(\omega,\varphi)=\frac{e^{\la\omega,\Psi_\lambda(\varphi)\ra}}{L_\lambda(\Psi_\lambda(
\varphi))}.\end{equation}  Here\rom,
\begin{align*}\Psi_\lambda(\varphi)&{:=}\frac{1}{\alpha-\beta}\,\log\bigg(\frac{1-\beta\varphi}{1-\alpha\varphi}\bigg),
\qquad \lambda\ne2,\\
\Psi_2(\varphi)&{:=}\frac{\varphi}{\varphi+1},\end{align*} is a
holomorphic ${\cal D}_\C$-valued function defined in a
neighborhood of zero in ${\cal D}_\C$ which is invertible and
satisfies $\Psi_\lambda(0)=0$\rom.\end{corollary}

\begin{remark}\label{hjdtr}\rom{ Corollary~\ref{hxxd} shows that
the system of polynomials $\la\wick n,f^{(n)}\ra$,
$f^{(n)}\in{\cal D }_\C^{\hotimes n}$, $n\in\Z_+$, is a
generalized Appell system in terms of \cite{KSS}, see also
\cite{ADKS,D,Kach,KSWY} and the references therein. We also refer
to \cite{bere2} and the references therein for the study of the
Appell systems via the  theory of generalized translation
operators.}

\end{remark}

\begin{remark}\label{zuaesdftz6}\rom{ Note that \begin{align*}
\Psi^{-1}_\lambda(\varphi)&{:=}\frac{e^{\alpha\varphi}-e^{\beta\varphi}}{\alpha
e^{\alpha\varphi}-\beta e^{\beta\varphi}}, \qquad \lambda\ne2,\\
\Psi^{-1}_2(\varphi)&{:=}\frac{\varphi}{1-\varphi}.\end{align*}
}\end{remark}

\begin{remark}\label{zgu3w4}\rom{In the one-dimensional case (Remarks~\ref{bhgfzt} and
\ref{fzfrzdzrz}), the function $$z\mapsto{\cal
L}_\lambda(z){:=}\int_\R e^{zx}\,{\frak m}_\lambda(dx)$$ is well
defined and holomorphic on a neighborhood of zero in $\C$ (for the
explicit formula, replace $iu$ in formula \eqref{tzrdt} with $z$).
Furthermore, for all $x\in\R$ and $z$ from the neighborhood of
zero, we have \cite{Meixner} $$ {\cal
G}_\lambda(x,z)=\frac{e^{x\Psi_\lambda(z)}}{{\cal
L}_\lambda(\Psi_\lambda (z))}, $$ where
$\Psi_\lambda(z)=\frac{1}{\alpha-\beta}\log\left(\frac{1-\beta z
}{1-\alpha z}\right)$, $\lambda\ne2$, $\Psi_2(z)=\frac z{z+1}$.

}\end{remark}

\noindent {\it Proof of Corollary~\rom{\ref{hxxd}}}. Fix
$\lambda\ge0$. Let $\widetilde L_\lambda(\varphi)$ denote the
right hand side of \eqref{jasuetfzcw}, respectively \eqref{hisav},
with $i\varphi$ replaced with $\varphi$.  It follows from the
proof of Theorem~\ref{z8fzvear} and \cite[sect.~3.2]{BK} that
there exists $ \tau_0\in T$ such that $\mu_\lambda({\cal
H}_{-\tau_0})=1$. By Proposition~\ref{ihaugcuzd}, there exists a
neighborhood
 of zero in ${\cal D}_\C$, denoted by $O_{\tau_0}$, such that for all
$\varphi\in O_{\tau_0}$  and all $\omega\in{\cal H}_{-\tau_0}$
\begin{equation}\label{jhvfdy}
e^{\la\omega,\varphi\ra}=\tilde L
_\lambda(\varphi)\,G_\lambda(\omega,\Psi^{-1}_\lambda(\varphi)).\end{equation}
Since   the number of the summands in the sum on the right hand
side of \eqref{1.5} is $n!$, we easily conclude that
\begin{equation}\label{uiafrsd}
\|G(\cdot,\Psi_\lambda^{-1}(\varphi))\|^2_{L^2({\cal
D}';\mu_\lambda)} =\sum_{n=0}^\infty \bigg(\frac 1{n!}\bigg)^2
\|(\Psi_\lambda^{-1}(\varphi))^{\otimes n}\|^2_{\Fext^{(n)}({\cal
H })}\,n!<\infty\end{equation} for all $\varphi$ from some (other)
neighborhood of zero in ${\cal D}_\C$, denoted by $O$. Then, for
all $\varphi\in \widetilde O_{\tau_0}{:=} O_{\tau_0}\cap O$, we
get by \eqref{jhvfdy} and \eqref{uiafrsd}:
\begin{align} \int_{{\cal D}'}e^{\la\omega,\varphi\ra}\,\mu_\lambda(d\omega)
&=\int_{{\cal
H}_{-\tau}}e^{\la\omega,\varphi\ra}\,\mu_\lambda(d\omega)\notag\\
&=\widetilde L_\lambda(\varphi)\int_{{\cal H}_{-\tau}}
G_\lambda(\omega,\Psi_\lambda^{-1}(\varphi))\,\mu_\lambda(d\omega)\notag\\
&=\widetilde L_\lambda(\varphi). \label{jhcs}\end{align} Formula
\eqref{ihugayzc} in the case of ${\cal H}_{-\tau_0}$ follows from
\eqref{jhvfdy} and \eqref{jhcs} The general case now easily
follows from Proposition~\ref{ihaugcuzd}.\quad $\blacksquare$

\section{Operators $\di_x$ and $\di^\dag_x$}\label{gtzdfzrdr}
For each $\tau \in T$, we introduce on ${\cal P}({\cal D}')$ a
Hilbertian norm $\|\cdot\|_{\lambda,\, \tau}$ as follows: for any
$\phi\in{\cal P}({\cal D}')$ of the form
$\phi(\omega)=\sum_{n=0}^N \la\wick n,f^{(n)}\ra$ (cf.\
Corollary~\ref{uiwef}), we set $$
\|\phi\|_{\lambda,\,\tau}^2{:=}\sum_{n=0}^N \|f^{(n)}\|_\tau^2\,
n!. $$ Let $({\cal D}_\lambda)_{\tau}$ denote the Hilbert space
obtained by closing ${\cal P}({\cal D}')$ in this  norm. By
\cite[Theorem~34]{KSS}, there exists $\tau_0\in T$ such that the
Hilbert space $({\cal D}_\lambda)_{\tau_0} $ is densely and
continuously embedded into $L^2({\cal D}';\mu_\lambda)$.  Just as
in Section~\ref{uisdygfuzuzg}, we first set $T'{:=}\{\tau\in
T:\tau\ge\tau_0\}$ and then re-denote $T{:=}T'$. Thus, each
$({\cal D }_\lambda)_\tau$, $\tau\in T$, consists of
($\mu_\lambda$-classes of) functions on ${\cal D}'$ of the form $$
\phi(\omega)=\sum_{n=0}^\infty \la \wick n,f^{(n)} \ra $$with
$f^{(n)}\in{\cal H}_\tau^{\hotimes n}$ and
\begin{equation*}\|\phi\|_{\lambda,\,\tau}^2{:=}\sum_{n=0}^\infty
\|f^{(n)}\|_\tau^{2}(n!)^2 <\infty. \end{equation*} Let
\begin{equation}\label{sjkbfb} ({\cal D}_\lambda){:=}\projlim
_{\tau\in T}({\cal D }_\lambda)_\tau,\end{equation} which is a
nuclear space \cite[Theorem~32]{KSS}. (We note that, though
 only the case of a nuclear space that is the projective
limit of a countable family of Hilbert space is considered in
\cite{KSS}, all the results we cite from this paper admit a
straightforward generalization to the case of a general nuclear
space.)

Denote by ${\cal E}_{\mathrm min}^1({\cal D}_\C')$ the set of all
entire functions on ${\cal D}_\C'$  of first order of growth and
of minimal type, i.e., a function $\phi$ entire on ${\cal D}_\C'$
belongs to ${\cal E}_{\mathrm min}^1({\cal D}_\C')$ if and only if
$$\forall \tau\in T,\ \forall\varepsilon >0\ \exists C>0:\ \forall
\omega\in{\cal H}_{-\tau,\,\C}:\quad |\phi(\omega)|\le C
e^{\varepsilon \|\omega\|_{-\tau}}.$$ Denote by ${\cal E}_{\mathrm
min}^1({\cal D}')$ the set of restrictions to ${\cal D}'$ of
functions from ${\cal E}_{\mathrm min}^1({\cal D}_\C')$. Following
\cite{KSS, KSWY}, we then  introduce  norms on ${\cal E}_{\mathrm
min}^1({\cal D}_\C')$, and hence on  ${\cal E}_{\mathrm
min}^1({\cal D}')$, as follows. For  each $\phi\in {\cal
E}_{\mathrm min}^1({\cal D}_\C')$ and for any $\tau\in T$ and
$q\in\N$, we set $$ n_{\tau,\,q}(\phi){:=}\sup_{z\in{\cal
H}_{-\tau,\,\C}}\big( |\phi(z)|\exp(-2^{-q}\,\|z\|_{-\tau}\big).$$
Next, each $\phi\in {\cal E}_{\mathrm min}^1({\cal D}_\C')$ can be
uniquely represented in the form $ \phi(z)=\sum_{n=0}^\infty \la
z^{\otimes n},f^{(n)}\ra$, and we set, for any $\tau\in T$ and
$q\in\N$, $$ N_{ \tau,\,q} (\phi){:=}\sum_{n=0}^\infty
\|f^{(n)}\|_\tau^2\,(n!)^2\, 2^{nq}.$$ By \cite[Theorems~2.5, 3.8
and subsec.~6.2]{KSS}, the three systems of norms on $({\cal D
}_\lambda)$: $$(\|\cdot\|_{\lambda,\,\tau},\  \tau\in T),\quad
(n_{\tau,\,q}(\cdot),\ \tau\in T,\, q\in\N),\quad
(N_{\tau,\,q}(\cdot),\ \tau\in T,\, q\in\N) ,$$ are equivalent,
and hence determine the same topology on $({\cal D }_\lambda)$.

As easily seen, for each $\varphi\in{\cal D}$, $I_\lambda
a_1^-(\xi)I_\lambda^{-1}$ can be extended to a continuous operator
on $({\cal D}_\lambda)$. We denote this operator $A_1^-(\xi)$.
Next, for each $x\in X$, we denote by $\di_x$ the linear
continuous operator on $({\cal D}_\lambda)$ defined by $$
\di_x\la\wick
n,f^{(n)}\ra{:=}n\la\wick{(n-1)},f^{(n)}(x,\cdot)\ra,\qquad \qquad
f^{(n)}\in{\cal D}_\C^{\hotimes n}. $$

\begin{lemma}\label{jkyxdvc} For any $\phi\in({\cal D}_\lambda)$\rom, we
have
 $$
\forall \omega\in{\cal D}':\qquad (A_1^-(\xi) \phi)(\omega)=\int_X
 \xi(x)(\di_x \phi)(\omega)\, \sigma(dx).$$
\end{lemma}

\noindent {\it Proof}. Let $\tau_0\in T$ be such that $\|\delta_x
\|_{\tau_0}\le1$ for all $x\in X$. Fix $\omega\in{\cal
H}_{-\tau}$, $\tau\ge\tau_0$. Let $\tau'>\tau$  be such that  the
inclusion ${\cal H}_{\tau_0}\hookrightarrow{\cal H}_\tau$  is
quasi-nuclear. By \cite[Proposition~22]{KSS}, we have for any
$\varepsilon>0$
\begin{equation}\label{hjsdfv}\|\wick n\|_{-\tau'}\le n!\,
C_\tau^n \exp(\varepsilon\|\omega\|_{-\tau}),\qquad
C_\tau>0.\end{equation} Hence, we may estimate $$
\sum_{n=1}^\infty n|\la\wick {(n-1)} ,f^{(n)}(x,\cdot)\ra|\le
\sup_{x\in X}\|\delta_x\|_{-\tau'}
\exp(\varepsilon\|\omega\|_{-\tau})\sum_{n=1}^\infty (n-1)!\,
C_\tau^{n-1} n\|f^{(n)}\|_{\tau'}. $$ Therefore, for each
$\phi=\sum_{n=0}^\infty \la{:}\cdot^{\otimes
n}{:}_\lambda,f^{(n)}\ra\in ({\cal D}_\lambda) $, we get, by the
majorized convergence theorem, $$\int_X \xi(x)\sum_{n=0}^\infty
n\la\wick {(n-1)}, f^{(n)}(x,\cdot)\ra\,\sigma(dx)=
\sum_{n=0}^\infty \la\wick {(n-1)},n\la f^{(n)},\xi\ra\,\ra,$$
where $$\la f^{(n)},\xi\ra(x_1,\dots,x_{n-1}){:=}\int_X
f^{(n)}(x,x_1,\dots,x_{n-1})\xi(x)\,\sigma(dx).$$ From here, the
lemma follows.  \quad $\blacksquare$

\begin{remark}\label{uisdf}\rom{Let $\tau,\tau'\in T$ be as in
proof of Lemma~\ref{jkyxdvc}. Note that the operator $\di_x$ acts
continuously in $({\cal D}_\lambda)_\tau$.  By \eqref{hjsdfv}, we
then have for all $\omega\in{\cal H}_{-\tau}$ and all
$\varphi\in{\cal D}_\C$ such that
$\|\varphi\|_{\tau'}<\min(1,C_\tau^{-1})$:
$$\di_x\sum_{n=0}^\infty \frac1{n!}\, \la\wick n,\varphi^{\otimes
n}\ra=\varphi(x)\sum_{n=0}^\infty \frac1{n!}\, \la\wick
n,\varphi^{\otimes n}\ra.$$ By Corollary~\ref{hxxd}, we get for
all $\omega\in{\cal H}_{-\tau}$ and all $\varphi$ from a
neighborhood of zero in ${\cal D}_\C$: $$ \di_x
e^{\la\omega,\Psi_\lambda(\varphi)\ra}=\varphi(x)
e^{\la\omega,\Psi_\lambda(\varphi)\ra},$$ and consequently
\begin{equation}\di_x
e^{\la\omega,\varphi\ra}=(\Psi_\lambda^{-1}(\varphi))(x)e^{\la\omega,\varphi\ra}.\label{jkyd}\end{equation}
Let $\nabla_x$ denote the G\^ateaux derivative of a function
defined on ${\cal D}'$ in direction $\delta_x$, i.e., $\nabla_x
F(\omega){:=}\frac d{dt}\big|_{t=0}F(\omega+t\delta_x)$. Clearly,
\begin{equation}\label{uztzewef}
\nabla_x
e^{\la\omega,\varphi\ra}=\varphi(x)e^{\la\omega,\varphi\ra}.\end{equation}
Comparing \eqref{jkyd} and  \eqref{uztzewef}, we get (at least
informally):
\begin{equation}\label{dtrdrt}\di_x=\Psi_\lambda^{-1}(\nabla_x).\end{equation}
}\end{remark}

\begin{remark}\label{aw<qw}\rom{ In the one-dimendsional case
(Remarks~\ref{bhgfzt}, \ref{fzfrzdzrz}, \ref{zgu3w4}), we define a
linear operator $A$ by $$ A
P^{(n)}_\lambda(\cdot)=nP_\lambda^{(n-1)}(\cdot).$$ By
Remarks~\ref{fzfrzdzrz} and \ref{zgu3w4}, one then gets (cf.\
\cite{Meixner})
\begin{equation}A=\Psi_\lambda^{-1}(D),\label{tdtrdst}\end{equation} where $D=\frac
d{dx}$. Thus, \eqref{dtrdrt} is an  infinite-dimensional
counterpart of \eqref{tdtrdst}.

}\end{remark}

\begin{theorem}\label{hjvcafh} For each $\lambda\ge0$ and for all
$\phi\in({\cal D}_\lambda)$ and $\omega\in{\cal D}'$\rom:
\begin{align}(\di_x \phi)(\omega)&=\int_\R
\big(\phi(\omega+s\delta_x)-\phi(\omega)\big)s\nu_\lambda(ds),\label{jhlxcvcg}\\
(A_1^-(\xi)\phi)(\omega)&=\int_{X\times\R}(\phi(\omega+s\delta_x)-\phi(\omega))s\xi(x)\,
\sigma(dx)\, \nu_\lambda(ds),\notag
\end{align}where $x\in X$\rom, $\xi\in{\cal D}$\rom, and $\nu_\lambda$ is the L\'evy measure of
$\mu_\lambda$  \rom(see Corollaries~\rom{\ref{audstz}} and
\rom{\ref{uiawfruza}}\rom)\rom.
\end{theorem}

\noindent{\it Proof}. By Lemma~\ref{jkyxdvc}, it suffices to prove
the statement only for $\di_x$. Fix $\tau\ge\tau_0$ as in proof of
Lemma~\ref{jkyxdvc}. Then, for all $\omega\in{\cal H}_{-\tau}$ and
all $\varphi$ from a neighborhood of zero in ${\cal D}_\C$, we
have
\begin{align} \int_\R \big(
e^{\la\omega+s\delta_x,\varphi\ra}-e^{\la \omega,\varphi\ra}
\big)s\, \nu_\lambda(ds)&=
e^{\la\omega,\varphi\ra}\int_\R\big(e^{s\varphi(x)}-1\big)s\nu_\lambda(ds)\notag\\
&=
e^{\la\omega,\varphi\ra}(\Psi_\lambda^{-1}(\varphi))(x),\label{hjacf}\end{align}
where $\Psi_\lambda^{-1}$ is given by \eqref{zuaesdftz6}. The
latter equality in \eqref{hjacf} can be derived by differentiating
in $\theta$ the following equality: \begin{align*}\int_\R
(e^{s\theta}-1-s\theta)\,\nu_\lambda(ds)&= \int_\R
\sum_{n=2}^\infty \frac{s^{n-2}\theta^n}{n!}\, s^
2\,\nu_\lambda(ds)\\
&=\begin{cases}\displaystyle-\frac1{\alpha\beta}\,\log\bigg(
\frac{\alpha e^{-\beta\theta}-\beta
e^{-\alpha\theta}}{\alpha-\beta}\bigg),& \lambda\ne2,\\
-\log(1-\theta)-\theta,&\lambda=2,\end{cases}\end{align*} which
holds for all $\theta$ from a neighborhood of zero in $\C$ (this
equality has been already used in course of the proof of
Corollaries~\ref{audstz} and \ref{uiawfruza}). Therefore, by
Remark~\ref{uisdf}, we have \eqref{jhlxcvcg} for all $x\in X $,
$\omega\in{\cal H}_{-\tau}$ and $\phi=G_\lambda(\cdot,\varphi)$,
where $\varphi$ runs through a neighborhood of zero in ${\cal
D}_\C$, denoted by $U_\tau$.

As easily seen, there exists $\varepsilon >0$ such that
\begin{equation}\label{hdask} \int_\R
e^{\varepsilon|s|}s^2\,\nu_\lambda(ds)<\infty\end{equation} (in
the case $\lambda<2$, see e.g.\ \cite[p.~180]{Chihara}).
 Choose $q\in\N$ such that  $2^{-q/2}<\varepsilon$. Choose
$\tau_1,\tau_2\in T$ and $q_1\in\N$ such that
$\tau<\tau_1<\tau_2$,
\begin{equation}\label{kjadg} C N_{\tau,\,q}(\cdot)\le C_1
n_{\tau_1,\,q_1}(\cdot)\le \|\cdot\|_{\lambda,\,\tau_2},\qquad
C,C_1>0.\end{equation}  Fix any $\phi\in({\cal D}_\lambda)$. For
each $k\in\N$,  let $\phi_k$ be a linear combination of functions
$G_\lambda(\cdot,\varphi)$ with $\varphi\in U_\tau$ and let
$\|\phi_k-\phi\|_{\lambda,\,\tau_2}\to0$ as $k\to\infty$
(evidently, such a sequence $\{\phi_k\}$ always exists). By
\eqref{kjadg}, we then have
\begin{equation}\label{ahdfgtgg}
n_{\tau_1,\,q_1}(\phi_k-\phi)\to0\quad\text{as
}k\to\infty.\end{equation} For a fixed $\omega\in{\cal
H}_{-\tau}\subset{\cal H}_{-\tau_1}$, we define $$\widetilde
\phi_k(z){:=}\phi_k(\omega+z)-\phi_k(\omega), \qquad \widetilde
\phi(z){:=}\phi(\omega+z)-\phi(\omega)$$ for $z\in{\cal
H}_{-\tau_1,\,\C}$. One easily checks that \eqref{ahdfgtgg}
implies $$ n_{\tau_1,\,q_1}(\widetilde \phi_k-\widetilde
\phi)\to0\quad\text{as }k\to\infty,$$ and hence, by \eqref {kjadg}
$$ N_{\tau,\,q}(\widetilde \phi_k-\widetilde
\phi)\to0\quad\text{as }k\to\infty.$$ We have $$\widetilde
\phi_k(z)=\sum_{n=1}^\infty \la z^{\otimes n},f^{(n)}_k\ra$$ (note
that $\widetilde \phi_k(0)=0$). Since $\|\delta_x\|_{-\tau}\le
\|\delta_x\|_{-\tau_0}\le1$ for all $x\in X$, we get, using the
Cauchy inequlity (compare with \cite[proof of Lemma~2.7]{KSWY})
\begin{align*} |\widetilde \phi_k(s\delta_x)|\, |s|^{-1}&\le
\sum_{n=1}^\infty |s|^{n-1}\|f_k^{(n)}\|_{\tau_1}\\ &\le
\sum_{n=1}^\infty {\max(1,|s|)}^n \|f_k^{(n)}\|_{\tau_1}\\ &\le
N_{\tau,\, q}(\widetilde\phi_k)\exp\big(2^{-q/2}\max(1,|s|)\big)\\
&\le C \exp\big(\varepsilon\max(1,|s|)\big),\qquad
k\in\N,\end{align*} where $C\in(0,\infty)$ is independent of
$k\in\N$. Hence, by \eqref{hdask} and the majorized convergence
theorem,
\begin{equation}\label{hgfflklk} \int_\R \widetilde
\phi_k(s\delta_x)\, \nu_\lambda(ds)\to \int_\R \widetilde
\phi(s\delta_x)\,\nu_\lambda(ds)\quad\text{as
}k\to\infty.\end{equation} Finally, $\di_x$ acts continuously on
$({\cal D}_\lambda)_{\tau_2}$, and hence  $\|\di_x \phi_k-\di_x
\phi\|_{\lambda,\,\tau_2}\to0$. Thus,  by \eqref{kjadg}, we get
$n_{\tau_1,\,q_1}(\di_x \phi_k-\di_x \phi)\to0$ as $k\to\infty$.
Therefore, $(\di_x \phi_k)(\omega)\to (\di_x \phi)(\omega)$ as
$k\to\infty$, which, together with \eqref{hgfflklk} concludes the
proof.\quad $\blacksquare$\vspace{2mm}

Let $({\cal D}_\lambda)_{-\tau}$ denote the dual space of $({\cal
D }_\lambda)_\tau$. Analogously to \eqref{hxvghyf}, we realize
$({\cal D}_\lambda)_{-\tau}$ as the Hilbert space consisting of
sequences $F=(F^{(n)})_{n=0}^\infty$, $F^{(n)}\in{\cal
H}_{-\tau,\,\C}^{\hotimes n}$, such that
$$\|F\|_{\lambda,\,-\tau}^2{:=}\sum_{n=0}^\infty\|F^{(n)}\|_{-\tau}^2<\infty,$$
with the scalar  product in $({\cal D}_\lambda)_{-\tau}$ generated
by the Hilbertian norm $\|\cdot\|_{\lambda,\,-\tau}$, and the dual
pairing of $F=(F^{(n)})_{n=0}^\infty$ with an element
\begin{equation}\label{svhgu}\phi=\sum_{n=0}^\infty \la{:}\cdot^{\otimes
n}{:}_\lambda,f^{(n)}\ra\in ({\cal
D}_\lambda)_{\tau}\end{equation} is given by $$ \la\!\la
F,\phi\ra\!\ra=\sum_{n=0}^\infty \la \overline{F^{(n)}},f^{(n)}
\ra\,n!.$$ By \eqref{sjkbfb}, we get the following representation
of the dual space of $({\cal D}_\lambda)$: $$ ({\cal
D}_\lambda)^*=\indlim_{\tau\in T}({\cal D}_\lambda)_{-\tau}.$$
Each test space $({\cal D}_\lambda)_\tau$ can be embedded into
$({\cal D}_\lambda)_{-\tau}$ by setting, for $\phi$ of the form
\eqref{svhgu}, $\iota(\phi){:=}(f^{(n)})_{n=0}^\infty$. In what
follows, we will just write $\phi$ instead of $\iota(\phi)$ to
simplify notations.

For each $x\in X$, we define an operator $\di_x^\dag$ on $({\cal D
}_\lambda)^*$ by $$\di_x^\dag
(F^{(n)})_{n=0}^\infty{:=}(\delta_x\hotimes F^{(n-1)})_{n=0}
^\infty.$$ Evidently, $\di_x^\dag$ is the dual operator of $\di_x$
and acts continuously in each $({\cal D}_{\lambda})_{-\tau}$,
$\tau\ge\tau_0$ (with $\tau_0$ as in the proof of
Lemma~\ref{jkyxdvc}).  We define operators
${:}\omega(x){:}_\lambda\cdot$ and $\omega(x)\cdot$, acting
continuously from each $({\cal D}_\lambda)_{\tau}$ into $({\cal
D}_\lambda)_{-\tau}$, $\tau\ge\tau_0$, by
\begin{equation}\label{uzftz}{:}\omega(x){:}_\lambda\cdot {:=}\di_x^\dag+\lambda
\di_x^\dag\di_x+\di_x+\di_x^\dag\di_x\di_x,\qquad
\omega(x)\cdot{:=}{:}\omega(x){:}_\lambda\cdot
+c_\lambda\operatorname{id}.\end{equation} Analogously to
$A^-_1(\xi)$, we define operators $A_2^-(\xi)$, $A^0(\xi)$, and
$A^+(\xi)$. Let also $\la\wick 1,\xi\ra\cdot$ and
$\la\omega,\xi\ra\cdot$ denote the operators of multiplication by
$\la\wick 1,\xi\ra$ and $\la\omega,\xi\ra$, respectively. We then
easily get the following integral representation:
\begin{gather}A^+(\xi)=\int_X
\sigma(dx)\,\xi(x)\di_x^\dag,\quad A^0(\xi)=\int_X
\sigma(dx)\,\xi(x)\di_x^\dag \di_x,\quad
  A_2^-(\xi) =\int_X
\sigma(dx)\,\xi(x)\di_x^\dag \di_x\di_x,\notag\\ \la\wick
1,\xi\ra\cdot= \int_X \sigma(dx)\,\xi(x)
{:}\omega(x){:}_\lambda\cdot,\qquad \la\omega,\xi\ra\cdot = \int_X
\sigma(dx)\,\xi(x) \omega(x)\cdot,\label{huadcf}
\end{gather} where the integrals are understood in the sense that
one applies pointwisely the integrand operator to a test function
from $({\cal D}_\lambda)$, then dualizes the result with another
test function, and finally integrates the obtained function with
respect to the measure $\sigma$.

Before formulating the next theorem, we note that, for each $x\in
X $, the operator $\nabla_x$ introduced  in Remark~\ref{uisdf}
acts continuously on $({\cal D}_\lambda)$. This can be easily
shown  using methods as in the proof of Theorem~\ref{hjvcafh}.
Furthermore, we introduce on $({\cal D}_\lambda)$ continuous
operators $\nabla_{\alpha-\beta,\,x}$ and ${\cal
U}_{\alpha-\beta,\,x}$ as follows:
\begin{align*}\nabla_{\alpha-\beta,\,x}
\phi(\omega)&{:=}\frac{\phi(\omega+(\alpha-\beta)\delta_x)-\phi(\omega)}
{\alpha-\beta},\\ {\cal
U}_{\alpha-\beta,\,x}\phi(\omega)&{:=}\phi(\omega-(\alpha-\beta)\delta_x).\end{align*}
Notice that, in the case $\lambda\in[0,2)$, we have
$\alpha-\beta=\alpha-\bar\alpha$, which is a purely imaginary
number, so when writing either $\phi(
\omega+(\alpha-\beta)\delta_x)$ or $\phi(
\omega-(\alpha-\beta)\delta_x)$, we understand under $\phi$ its
entire extension to $ {\cal D}_\C'$.

\begin{theorem}\label{jhlsdg} For each $x\in X$\rom, $\di^\dag_x$ considered as an operator
from $({\cal D}_\lambda)$ into $({\cal D}_\lambda)^*$ has the
following representation\rom:
\begin{equation}\label{sdouh} \di_x^\dag=\begin{cases}
{:}\omega(x){:}_\lambda\cdot
(1+\alpha\nabla_{\alpha-\beta,\,x})^2{\cal
U}_{\alpha-\beta,\,x}-(1+\alpha\nabla_{\alpha-\beta,\,x})\nabla_{\alpha-\beta,\,x}\,{\cal
U }_{\alpha-\beta,\,x},&\lambda\ne2,\\ \omega(x)\cdot
(\nabla_x-1)^2+(\nabla_x-1),&\lambda=2.\end{cases}\end{equation}
\end{theorem}

\noindent {\it Proof}. We prove the theorem only in the case
$\lambda\ne2$
 and refer to \cite[Lemma~7.1]{KL} for the case  $\lambda=2$. By
 \eqref{uzftz} and \eqref{jkxhgvfuz},
\begin{align} {:}\omega(x){:}_\lambda\cdot
&=\di_x^\dag(1+\lambda\di_x+\di_x^2)+\di_x\notag\\
&=\di^\dag_x(1-\alpha\di_x)(1-\beta\di_x)+\di_x.\label{hjdsegeg}\end{align}
Fix $\tau\ge\tau_0$ and choose $\tau'\in T$, $\tau'>\tau$, such
that the operators $(1+\alpha\nabla_{\alpha-\beta,\,x})^2{\cal
U}_{\alpha-\beta,\,x}$ and
$(1+\alpha\nabla_{\alpha-\beta,\,x})\nabla_{\alpha-\beta,\,x}\,{\cal
U }_{\alpha-\beta,\,x}$ act continuously from $({\cal
D}_\lambda)_{\tau'}$ into $({\cal D}_\lambda)_{\tau}$. Take any
$\varphi\in{\cal D}$ such that $e^{\la\cdot,\varphi\ra}\in({\cal D
}_\lambda)_{\tau'}$.  Then, by \eqref{jkyd} and \eqref{hjdsegeg},
\begin{equation}\label{hjvsdvgvgs} {:}\omega(x){:}_\lambda\cdot
e^{\la\cdot,\varphi\ra}=\di_x^\dag
\big(1-\alpha(\Psi_\lambda^{-1}(\varphi))(x)\big)\big(1-\beta(\Psi
_\lambda^{-1}(\varphi))(x)\big)e^{\la\cdot,\varphi\ra}+
(\Psi_\lambda^{-1}(\varphi))(x)e^{\la\cdot,\varphi\ra}.
\end{equation}
Denoting
$$\nabla_{\alpha-\beta}\varphi(x){:=}\frac{e^{\varphi(x)(\alpha-\beta)}-1}{\alpha-\beta},$$
we get
\begin{equation}\label{hgasdghsd}(\Psi_\lambda^{-1}(\varphi))(x)
=\frac{\nabla_{\alpha-\beta}\varphi(x)}{1+\alpha
\nabla_{\alpha-\beta}\varphi(x) } \end{equation} (compare with
\cite[formula~(7.2)]{Meixner}).  By \eqref{hjvsdvgvgs}  and
\eqref{hgasdghsd}, $$ {:}\omega(x){:}_\lambda\cdot
e^{\la\cdot,\varphi\ra}=\di_x^\dag
\,\frac{e^{\varphi(x)(\alpha-\beta)}}{(1+\alpha\nabla_{\alpha-\beta}\varphi(x))^2}\,
e^{\la\cdot,\varphi\ra} +\frac{\nabla_{\alpha-\beta}\varphi(x)}
{1+\alpha\nabla_{\alpha-\beta}\varphi(x)}\,e^{\la\cdot,\varphi\ra},
$$ which yields \begin{multline}\label{hbsdhhs}
{:}\omega(x){:}_\lambda\cdot
(1+\alpha\nabla_{\alpha-\beta}\varphi(x))^2e^{-\varphi(x)(\alpha-\beta)}
e^{\la\cdot,\varphi\ra}\\=\di_x^\dag e^{\la\cdot,\varphi\ra}
+\nabla_{\alpha-\beta}\varphi(x)(1+\alpha\nabla_{\alpha-\beta}\varphi(x))e^{-\varphi(x)(\alpha-\beta)}
e^{\la\cdot,\varphi\ra}.
\end{multline} Since the set of all functions
$e^{\la\cdot,\varphi\ra}$ with $\varphi$ as above is total in
$({\cal D}_\lambda)_{\tau'}$ and since
${:}\omega(x){:}_\lambda\cdot$ and $\di_x^\dag$ act continuously
from $({\cal D}_\lambda)_\tau$ into $({\cal D}_\lambda)_{-\tau}$,
\eqref{hbsdhhs} implies \eqref{sdouh}. \quad $\blacksquare$

\begin{corollary}\label{uhsduhisuhis} For any $\phi\in({\cal
D}_\lambda)$ and $\xi\in{\cal D}$\rom, we have for all
$\omega\in{\cal D}'$
\begin{multline}\label{haszgtafza}A^+(\xi)\phi(\omega)=\la
\omega(x)
-c_\lambda,\xi(x)(1+\alpha\nabla_{\alpha-\beta,\,x})^2{\cal U
}_{\alpha-\beta,\, x} \phi(\omega)\ra \\ \text{}-\la
1,\xi(x)(1+\alpha\nabla_{\alpha-\beta,\,x})\nabla_{\alpha-\beta,\,x}\,{\cal
U }_{\alpha-\beta,\, x} \phi(\omega)\ra,\qquad\lambda \ne2,
\end{multline}
and $$A^+(\xi)\phi(\omega)=\la \omega(x)
,\xi(x)(\nabla_x-1)^2\phi(\omega)\ra+\la
1,\xi(x)(\nabla_x-1)\phi(\omega)\ra,\qquad \lambda=2,$$ where $x$
denotes the variable in which the dualization is carried out\rom.
\end{corollary}

\noindent{\it Proof}. Again, we prove only in the case
$\lambda\ne2$ and refer to \cite[Theorem~7.1]{KL} for the case
$\lambda=2$. Fix any  $\tau\ge\tau_0$. Using methods as in the
proof of Theorem~\ref{hjvcafh}, we conclude the existence of $
\tau_1\in T$, $\tau_1>\tau$, such that, for any fixed
$\omega\in{\cal H }_{-\tau}$ and  $z\in\C$,
$|z|\le2|\alpha-\beta|$, and for any sequence $\{\phi_k,\
k\in\N\}\subset ({\cal D}_\lambda)_{\tau_1}$ such that
$\phi_k\to\phi$ in  $({\cal D}_\lambda)_{\tau_1}$ as $k\to\infty$,
 we have $\xi\widetilde\phi_{k,\,z}\to \xi\widetilde\phi_z$ in ${\cal
H}_{\tau,\,\C}$ as $k\to\infty$. Here, $$X\ni x\mapsto
\widetilde\phi_{k,\,z}(x){:=}\phi_k(\omega+z\delta_x)\in\C,\qquad
X\ni x\mapsto
\widetilde\phi_z(x){:=}\phi(\omega+z\delta_x)\in\C.$$ Choose
$\tau_2>\tau_1$ such that $C
n_{\tau_1,1}(\cdot)\le\|\cdot\|_{\lambda,\,\tau_2}$, $C>0$, and
choose $\tau_3>\tau_2$ such that $A^+(\xi)$ acts continuously from
$({\cal D}_\lambda)_{\tau_3}$ into $({\cal
D}_{\lambda})_{\tau_2}$. Fix any $\varphi\in{\cal D}$ such that
$e^{\la\cdot,\varphi\ra}\in({\cal D}_\lambda)_{\tau_3}$. Then, for
any $\phi\in({\cal D}_\lambda)$, we get by (the proof of)
Theorem~\ref{jhlsdg} and \eqref{huadcf} \begin{gather*} \la\!\la A
^+(\xi)e^{\la\cdot,\varphi\ra},\psi \ra\!\ra=\int_X \xi(x )
\la\!\la \di_x^\dag e^{\la\cdot,\varphi\ra},\psi
\ra\!\ra\,\sigma(dx)\\=
\int_X\xi(x)(1+\alpha\nabla_{\alpha-\beta}\varphi(x))^2
e^{-\varphi(x)(\alpha-\beta)}\la\!\la {:}\omega(x){:}_\lambda\cdot
e^{\la\cdot,\varphi\ra},\psi\ra\!\ra\,\sigma(dx)
\\ \text{}- \int_X
\xi(x)(1+\alpha\nabla_{\alpha-\beta}\varphi(x))\nabla_{\alpha-\beta}\varphi(x)
e^{-\varphi(x)(\alpha-\beta)}\la\!\la
e^{\la\cdot,\varphi\ra},\psi\ra\!\ra\,\sigma(dx)
\\ =\la\!\la \,\la {:}\cdot^{\otimes 1}{:}_\lambda,\xi (1+\alpha\nabla_{\alpha-\beta}\varphi)^2
e^{-\varphi(\alpha-\beta)}\ra\cdot
e^{\la\cdot,\varphi\ra},\psi\ra\!\ra\,\sigma(dx) \\ \text{}-
\la\!\la\,\la
1,\xi(1+\alpha\nabla_{\alpha-\beta}\varphi)\nabla_{\alpha-\beta}\phi\,
e^ {-\varphi(\alpha-\beta)}\ra
e^{\la\cdot,\varphi\ra},\psi\ra\!\ra,
\end{gather*}
which implies \eqref{haszgtafza} for
$\phi=e^{\la\cdot,\varphi\ra}$ and $\omega\in{\cal H}_{-\tau}$.
Now, approximate an arbitrary $\phi\in({\cal D}_\lambda)$ in the
$\|\cdot\|_{\lambda,\,\tau_3}$ norm by a sequence
$\{\phi_k,\,k\in\N\}$ of linear combinations of exponents as
above.  Then, $A^+(\xi)\phi_k\to A^+(\xi)\phi$ as $k\to\infty$ in
the $\|\cdot\|_{\lambda,\,\tau_2}$ norm, and hence in the
$n_{\tau_1,\,1}(\cdot)$ norm. In particular,
$A^+(\xi)\phi_k(\omega)\to A^+(\xi)\phi(\omega)$ for any
$\omega\in{\cal H}_{-\tau_1}$. Furthermore, for any fixed
$\omega\in{\cal H}_{-\tau}$, we get \begin{align*}
\xi(1+\alpha\nabla_{\alpha-\beta,\bullet})^2{\cal
U}_{\alpha-\beta,\bullet}\phi_k(\omega)&\to
\xi(1+\alpha\nabla_{\alpha-\beta,\bullet})^2{\cal
U}_{\alpha-\beta,\bullet}\phi(\omega),\\ \xi(
1+\alpha\nabla_{\alpha-\beta,\bullet})
\nabla_{\alpha-\beta,\bullet}{\cal
U}_{\alpha-\beta,\bullet}\phi_k(\omega)&\to  \xi(
1+\alpha\nabla_{\alpha-\beta,\bullet})
\nabla_{\alpha-\beta,\bullet}{\cal
U}_{\alpha-\beta,\bullet}\phi(\omega)
\end{align*}
in ${\cal H}_{-\tau,\,\C}$. From here, we evidently get
\eqref{haszgtafza} for an arbitrary $\phi\in({\cal D}_\lambda)$
and an arbitrary $\omega\in{\cal H}_{-\tau}$. \quad $\blacksquare$

\begin{remark}\rom{ Analogously to Corollary~\ref{uhsduhisuhis},
one can derive from Theorems~\ref{hjvcafh}, \ref{jhlsdg} and
\eqref{huadcf} explicit formulas for the action of the operators
$A^0(\xi)$ and $A^-_2(\xi)$ (see also \cite[Theorems~7.2, 7.3]{KL}
for the case $\lambda=2$).}\end{remark}

\noindent {\bf Acknowledgements}

I am grateful to Yu.~Berezansky for  many remarks on a preliminary
version of the paper and making me aware of the recent paper
\cite{bere3}. I am indebted to Yu.~Kondratiev and M.~J.~Oliveira
for useful discussions. I would  also like to thank the referee of
the paper for the suggestions on improvement of the first version
of the paper. The financial support of the SFB 611 and the DFG
Research Projects 436 RUS 113/593 is gratefully acknowledged.

\end{document}